\documentclass[reqno]{amsart}
\usepackage{amsmath}
\usepackage{amssymb}
\usepackage{stmaryrd}
\title{
A Simple Regularization of Hypergraphs
}
\author{Yoshiyasu Ishigami}
\address{Department of Information and Communication Engineering,
The University of Electro-Communications, Chofu, Tokyo 182-8585, Japan.}
\def\picture #1 by #2 (#3){
                \vbox to #2{
                        \hrule width #1  height 0pt depth
0pt
                        \vfill
                        \special{picture #3}}}
\def\scaledpicture #1 by #2(#3 scaled #4){{
                \dimen0=#1 \dimen1=#2
                \divide\dimen0 by 1000 \multiply\dimen0 by
#4
                \divide\dimen1 by 1000 \multiply\dimen1 by
#4
                \picture\dimen0 by \dimen1 (#3 scaled #4)}}
\newtheorem{adf}{Definition}[section]

\newtheorem{thm}{Theorem}[section]
\newtheorem{cor}[thm]{Corollary}
\newtheorem{lemma}[thm]{Lemma}
\newtheorem{prop}[thm]{Proposition}

\newtheorem{fact}[thm]{Fact}

\newtheorem{remark0}[thm]{Remark}
\newenvironment{remark}{%
\begin{remark0}\rm}{\thqed\end{remark0}}
\newtheorem{algorithm0}[thm]{Algorithm}
\newtheorem{asett}[thm]{Setup}

\newenvironment{df}{
\begin{adf}\begin{sl}}
{\end{sl}\thqed\end{adf}}
\newenvironment{sett}{
\begin{asett}\begin{sl}}{\end{sl}\thqed\end{asett}}
\newtheorem{example0}[thm]{Example}
\newenvironment{example}{%
\begin{example0}\rm}{\thqed\end{example0}}

\newcommand{\thqed}{\hfill\fbox{}\\ }
\newcommand{\Proof }{{\bf Proof} : }

\newcommand{\factqed}{\hfill\bigskip\fbox{}\\}

\newcommand{\lemmaqed}{\hfill\bigskip\fbox{}\\}

\def\brkt#1{\left({#1}\right)}
\def\dbrkt#1{\llbracket{#1}\rrbracket}

\def\ang#1{\langle{#1}\rangle}

\newcommand{\realset}{{\mathbb R}}

\newcommand{\naturalset}{{\mathbb N}}
\newcommand{\integerset}{{\mathbb Z}}
\newcommand{\Prob}{{\mathbb P}}
\newcommand{\Ex}{{\mathbb E}}

\newtheorem{procedure0}{Procedure}

\newcount\timehh\newcount\timemm\timehh=\time\divide\timehh by 60
\timemm=\time\count255=\timehh\multiply\count255 by-60
\advance\timemm by \count255
\ifnum\timehh=12\def\apm{pm}\else
\ifnum\timehh>12\def\apm{pm}\advance\timehh by-12\else
\def\apm{am}\fi\fi
\def\timestamp{\number\timehh\,:\,\ifnum\timemm<10 0\fi\number\timemm\,\apm}
\begin{document}
\maketitle
\begin{abstract}
We give a simple and natural (probabilistic) construction of hypergraph
regularization. 
It is done just by 
taking a constant-bounded number of random vertex samplings only one time (thus, iteration-free).
It is independent from the definition of quasi-randomness and 
 yields a new elementary proof of a strong hypergraph regularity lemma. 
Consequently, as an example of its applications, we have a
new self-contained proof of Szemer\'edi's classic theorem on
arithmetic progressions (1975) as well as its multidimensional
extension by Furstenberg-Katznelson (1978).
\end{abstract}
\section{Introduction}
\subsection{Szemer\'edi-type density theorems}
The following is often considered as one of the deepest theorems
in combinatorics.
\begin{thm}[Multi-dimensional Szemer\'edi Theorem --
Furstenberg-Katznelson (1978)\cite{FK78}]\label{24} For any
$\delta>0, \,r\ge 1$\,, and $F\subset \integerset^r$ with
$|F|<\infty$, if an integer $N$ is sufficiently large then for any
subset $S\subset \{0,1,\cdots,N-1\}^r$ with $|S|\ge \delta N^r$
there exist $a\in \{0,1,\cdots,N-1\}^r$ and $c\in [N]$ with
$a+cF\subset S$.
\end{thm}
Furstenberg and Katznelson (1978) \cite{FK78} proved this by using
ergodic theory. The special case of  $r=2$ and
$F=\{(0,0),(0,1),(1,0),(1,1)\}$ was first conjectured by
R.L.~Graham in 1970 (\cite{AS74,E73}). The case of $r=2$ and
$F=\{(0,0),(0,1),(1,0)\}$, was investigated initially by
Ajtai-Szemer\'edi (1974) \cite{AS74}.
\par
The following was first conjectured by Erd\H{o}s and Tur\'an (1936) \cite{ET36}.
\begin{cor}[Szemer\'edi (1975)\cite{Sz75}]\label{bnj08Sz}
For any $\delta>0$ and $m\ge 1$, there exists an integer $N$ such that
any subset $S\subset [N]$ with $|S|\ge \delta N$
contains an arithmetic progression of length $m.$
\end{cor}
Green and Tao \cite{GT06p} recently proved the existence of
arbitrarily long arithmetic progressions in the primes, in which
they used Szemer\'edi's theorem. 
\subsection{A brief history of hypergraph regularity}
Inspired by the success of the celebrated Graph Regularity Lemma
\cite{Sz}, research on quasi-random hypergraphs was initiated independently 
by at least four groups: 
Chung or Chung-Graham \cite{Chung,Chung91,CG,CG2,CG3},
Frankl-R\"{o}dl\cite{FR92}, 
Haviland-Thomason \cite{HT,HT2}, 
and Steger\cite{Steger90}(see \cite{PS92b} for its application).
For other earlier work, see
\cite{BNS,CT}. Also, Frankl-R\"{o}dl (2002) \cite{FR02} gives a
regularity lemma for $3$-uniform hypergraphs.
\par
Then R\"{o}dl and his collaborators \cite{RSk04,NRS}
and Gowers \cite{G} independently obtained their hypergraph
regularity lemmas. Slightly later, Tao \cite{Tao06} gave another
regularity lemma. 
\par
It has been noted that unlike the situation for graphs, there are
several ways one might define regularity for hypergraphs
(R\"{o}dl-Skokan \cite[pp.1]{RSk04},Tao-Vu \cite[pp.455]{TV}).
(For sparse hypergraphs, an essential difference appears. See
\cite[\S 10]{G3u}.) 
Kohayakawa et al. \cite[pp.188]{KNR} say that
the basic objects involved in the Regularity Lemma and the
Counting Lemma are already somewhat technical and that simplifying
these lemmas would be of great interest. In this paper we try to meet 
these requirements. We can naturally obtain strong quasi-random
properties not from one basic quasi-random property but from our
construction of a certain partition which we will define.
\par
In this paper, we give a new construction of hypergraph
regularization. Our regularization is achieved by a quite simple (probabilistic)
construction which makes it easy to understand why it works. Note
that our construction of regularization is new even if we assume 
we are working with ordinary graphs. 
In our construction, the number of random vertex samplings is not a fixed constant 
and our construction is iteration-free. 
(In later sections, we will see how different it is from property test more.)
But once the statement of our construction is given, its proof may be deduced naturally.
\par
For applications of the main result of this paper, see \cite{I06m,I06lr,I07}.
\subsection{Differences from the previous hypergraph regularities}
A Regularity lemma works well for applications when its 
counting lemma accompanies it. All of the previous proofs go as
follows.
\\
(i) Define regularity (a basic quasi-random property) for each
{\em cell} (a $k$-uniform $k$-partite hypergraph),
\\
(ii) Prove the existence of a partition in which most cells satisfy the regularity.
[Regularity Lemma]\\
(iii) 
Estimate the number of copies of a fixed 
hypergraph. [Counting Lemma]
\par
Our program will go as follows.
\\
(i') Define the construction of a partition. (Its existence will
be clear.)
\\
(ii') Estimate the number of copies of a fixed (colored) hypergraph.
\par
Once the definition of the construction via random samplings is given, 
the concept of our proof is simple.
The most interesting technical 
part in our proof is to use \lq linearity of expectation.\rq\ 
All of the
previous proofs use the dichotomy (or energy-increment)
explicitly and iteratively. (See \cite[\S6]{G},
 \cite[\S1]{Tao06a}.) Namely, when proving (ii), they define an \lq energy\rq\ (or index)
by the {\em supremum} (or maximum) of some (energy) function. (For
example, see \cite[eq.~(8)]{Tao06}.) It corresponds to
(\ref{061220}) in this paper. They consider the {\em supremum}
value of this energy over all subdivisions in each step. If the
energy significantly increases by {\em some} subdivision, they
take the {\em worst} subdivision as the base partition of the next
step. They then repeat this process. Since the energy is bounded,
this operation must stop at some step, in which case
 there is no quite bad subdivision, and thus, most cells should be quasi-random
(dichotomy).
\par
On the other hand, we (implicitly) take an {\em average}
subdivision instead of the worst one. The definition of our
regularization determines the probability space of partitions
(subdivisions). We also randomly decide on the number of vertex
samples to choose.
With these ideas, we can {\em hide} the troublesome dichotomy
iterations inside linear equations of expectations (\ref{lj20}). (
Imagine what would happen in (\ref{lj20}) if we replaced
$\Ex_\varphi$ by $\sup_\varphi$ in (\ref{061220}). ) (One of the
main reasons why Tao's \cite{Tao06} proof is relatively shorter
than the earlier two may be that he also reduced double-induction
concerns by preparing two partitions (coarse/fine) instead of one partition 
in each level $i\in [k-1].$ So in this sense, his regularity lemma
is seemingly weaker but still strong enough for proving 
removal lemmas and applications,  
which was his main interest. 
\par
We have two reasons why we will deal with multi-colored
hypergraphs instead of ordinary hypergraphs, even though almost
all previous researchers dealt with the usual 
hypergraphs (with black\&white edges).
 First, our proof of the regularity lemma will be natural.
Second, we can naturally combine subgraph (black\&invisible) and
induced-subgraph (black\&white) problems when we apply our result,
while the two have usually been discussed separately.
The set of these definitions to state our main theorem  
is new and helps us to simplify the arguments that follow. 
The magnitude of this effect is not small.
It is not hard for advanced readers to imagine that 
it would become even larger 
when we consider applications of our main theorem to other problems, 
some of which 
require to modify the proof of our main theorem itself.
\par
\allowdisplaybreaks[4]
\section{Statement of the Main Theorem}
In this paper, $\Prob$ and $\Ex$ will denote probability and
expectation, respectively. We denote conditional probability and
exepctation by $\Prob[\cdots|\cdots]$ and $\Ex[\cdots|\cdots].$
\begin{sett}\label{f1231}
Throughout this paper, we fix a positive integer $r$ and an \lq
index\rq\ set $\mathfrak{r}$ with $|\mathfrak{r}|=r.$ Also we fix
a probability space $({\bf \Omega}_i,{\mathcal B}_i,\Prob)$ for
each $i\in \mathfrak{r}$. We assume that ${\bf \Omega}_i$ is
finite (but its cardinality will not be a constant in our statements) 
and that ${\mathcal B}_i=2^{{\bf \Omega}_i}$ (for the sake of simplicity).
Write ${\bf \Omega}:=({\bf \Omega}_i)_{i\in \mathfrak{r}}$.
\end{sett}
In order to avoid using measure-theoretic jargon such as
measurability or Fubini's theorem, for the benefit readers who are
interested only in applications to discrete mathematics, we assume
${\bf \Omega}_i$ to be a (non-empty) finite set. However, our
arguments should be extendable to a general probability space. For
applications, ${\bf \Omega}_i$ usually would contain a huge number
of vertices, though we will not use this assumption in our proof.
(Note that this assumption has been actively used by many 
researchers.)
\par
For an integer $a$, we write $[a]:=\{1,2,\cdots,a\},$ and
${\mathfrak{r}\choose [a]}:=\dot{\bigcup}_{i\in
[a]}{\mathfrak{r}\choose i} =\dot{\bigcup}_{i\in [a]}\{I\subset
\mathfrak{r}| |I|=i\}.$ 
We also use the notation $[a,b]:=\{a,a+1,\cdots,b\}$ for 
integers $a,b.$
\begin{df}[(Colored hyper)graphs]\label{070724}
\quad
Suppose Setup \ref{f1231}. A {\bf $k$-bounded $(b_i)_{i\in
[k]}$-colored ($\mathfrak{r}$-partite hyper)graph} $H$ is an object with the following 
three ingredients 
:\\
$\bullet$ A union $V(H)=\dot{\bigcup}_{i\in \mathfrak{r}}V_i(H)$ of 
disjoint sets. The sets $V_i(H)$ and their elements are called  
{\bf vertex sets} and {\bf vertices} of $H$, respectively.
Write $V_J(H):=\{e\subset \dot{\bigcup}_{i\in J}V_i(H)\,: |e\cap
V_j(H)|=1\,, \forall j\in J\}$ whenever $J\subset {\mathfrak r}$.
Each element $e\in V_I(H)$ with $I\in {\mathfrak{r}\choose [k]}$ is called
 an {\bf (index-$I$ size-$|I|$) edge}.
\\
$\bullet$ For each $I\in {\mathfrak{r}\choose [k]}$, a set ${\rm C}_I(H)$ of 
exactly $b_{|I|}$ elements, where the elements are called 
 {\bf (face-)colors (of index $I$ and size $|I|$)}.
\\
$\bullet$ For each $I\in {\mathfrak{r}\choose [k]}$, a function from $V_I(H)$ to ${\rm 
C}_I(H)$. Denote by $H(e)$ the image of $e\in V_I(H)$ via the function.
\par
Let $I\in {\mathfrak{r}\choose [k]}$ 
and $e\in V_I(H).$ 
For another index $\emptyset\not=
J\subset I$, we denote by $e|_J$ the index-$J$ edge 
$e\setminus \brkt{\bigcup_{j\in I\setminus J}X_j}\in V_J(H)$.
We define the {\bf frame-color}
and {\bf total-color} of $e$ by vector ${H}(\partial
e):=({H}(e|_J)|\,\emptyset\not=J\subsetneq I)$ and by vector 
${H}(\ang{e})=H\ang{e}:=({H}(e|_J)|\,\emptyset\not=J\subset I).$
Write ${\rm TC}_I(H):=\{H\ang{e}|\,{e}\in X_I\},$ ${\rm
TC}_s(H):=\bigcup_{I\in {\mathfrak{r}\choose s}}{\rm TC}_I(H),$
and ${\rm TC}(H):=\bigcup_{s\in  [k]}{\rm TC}_s(H).$
\end{df}
\begin{example}
An ordinary ($\mathfrak{r}$-partite) graph is a $2$-bounded $(b_1,b_2)$-colored hypergraph 
with $b_1=1$ and $b_2=2$.
\par 
A triple $e=\{v_1,v_3,v_4\}$ of vertices is an index-$\{1,3,4\}$ edge 
if and only if $v_1\in X_1, v_3\in X_3$ and 
$v_4\in X_4.$ In any 
$k$-bounded $\mathfrak{r}$-partite hypergraph, 
any vertex in $X_i$ is an index-$\{i\}$ edge (whenever $k\ge 1$).
For two $k$-bounded $\mathfrak{r}$-partite hypergraphs $H$ and $H'$ 
with a common vertex set 
$V(H)=V(H')=\dot{\bigcup}_{i\in\mathfrak{r}}X_i$, 
all the edges of $H$ are also the edges of $H'$.
In this sense, our definition of the word \lq edge\rq\ 
is different from that in 
the classical (hyper)graph theory.
In our setting, the essential structure of a colored hypergraph 
is determined not by the set of edges but by the map 
from the edges to the colors. 
\par
All index-$I$ edges are colored not only when $|I|=k$ 
but also when $1\le |I|<k$, 
which is the reason why we call the hypergraph $k$-bounded 
instead of $k$-uniform.
\par
If $I=\{1,3,5\}, J=\{1,5\}, 
v_1\in X_1,v_3\in X_3, v_5\in X_5$ and 
$e=\{v_1,v_3,v_5\}$ then 
$e|_J=\{v_1,v_5\}.$ 
\end{example}
\par
Throughout the paper, 
we will try to embed an $r$-partite graph $S$ to another larger $r$-partite 
graph ${\bf G}$, where 
the $r$ vertex-sets of the larger graph will be always $({\bf \Omega}_i)_{i\in\mathfrak{r}}.$
And the larger graph and its vertices and edges will be denoted by bold fonts (ex. ${\bf G}, {\bf v}, 
{\bf v'}, {\bf e}, \cdots$) in order to avoid confusing them with those of the smaller graph.
The smaller graph will be always a simplicial-complex defined below.
\begin{df}[Simplicial-complexes]
A {\bf ($k$-bounded) simplicial-complex} is a $k$-bounded
(colored ${\mathfrak r}$-partite hyper)graph
such that for each $I\in {\mathfrak{r}\choose [k]}$
there exists at most one index-$I$ color called \lq invisible\rq\
and that if (the face-color of) an edge $e$ is invisible then
 (the face-color of) any edge $e^*\supset e$ is invisible. 
We call an edge invisible when the face-color of the edge is invisible.
An edge or its color is {\bf visible} if it is not invisible.
\par
For a $k$-bounded graph ${\bf G}$ on ${\bf \Omega}$ and $s\le k$,
let ${\mathcal S}_{s,h,{\bf G}}$ be the set of $s$-bounded
simplicial-complexes $S$ such that:
\\ (1) each of the $r$ vertex-sets of $S$ 
contains exactly $h$ vertices, and that, 
\\
(2) for 
$I\in
{\mathfrak{r}\choose [s]}$
 there is an injection from the index-$I$ 
visible colors of $S$ to the
index-$I$ colors of ${\bf G}$.\\ (When a visible color
$\mathfrak{c}$ of $S$ corresponds to another color $\mathfrak{c}'$
of ${\bf G}$, we simply write $\mathfrak{c}=\mathfrak{c}'$ without
presenting the injection explicitly.) For $S\in {\mathcal
S}_{s,h,{\bf G}}$, we denote by ${\mathbb V}_I(S)$ the set of
index-$I$ visible edges. Write ${\mathbb V}_i(S):=\bigcup_{I\in
{\mathfrak{r}\choose i}}{\mathbb V}_I(S)$ and ${\mathbb
V}(S):=\bigcup_i {\mathbb V}_i(S).$
\end{df}
For our purpose of this paper, all of the colors in the larger graph ${\bf G}$ can be 
considered to visible, though we will not use it logically.
\begin{df}[Partitionwise maps]
A {\bf partitionwise map} $\varphi$ is
a map
 from $r$ vertex sets $W_i,i\in \mathfrak{r},$ with $|W_i|<\infty$, to {\em the}
$r$ vertex sets (probability spaces) $ {\bf
\Omega}_i,i\in\mathfrak{r}$, such that each $w\in W_i$ is mapped
into ${\bf \Omega}_i.$ We denote by
$\Phi((W_i)_{i\in\mathfrak{r}})$ or
$\Phi(\bigcup_{i\in\mathfrak{r}}W_i)$ the set of partitionwise
maps from $(W_i)_i.$ When $W_i=\{(i,1),\cdots, (i,h)\}$ 
or  when $W_i$ are obvious and
$|W_i|=h$, we denote it by $\Phi(h)$. 
We write $
 \varphi({\mathbb D}):=
\dot{\bigcup}_{i\in \mathfrak{r}}\varphi(W_i)$ 
for $\varphi\in \Phi(\bigcup_{i\in\mathfrak{r}}W_i)$ 
(when we want to denote the range without saying the domain explicitly).
A partitionwise map is {\bf random} if and only if for every $i,$ each $w\in W_i$
is mutually independently mapped to a point in the probability
space ${\bf \Omega}_i$.
\par
Define $\Phi(m_1,\cdots,m_{k-1}):=\Phi(m_1)\times\cdots\times\Phi(m_{k-1}).$
\par
For two partitionwise maps $\phi\in \Phi((W_i)_i)$ and $\phi'\in \Phi((W'_i)_i),$ 
denote by $\phi\dot{\cup} \phi'$ 
the partitionwise map $\phi^*\in \Phi((W_i\dot{\cup} W'_i)_i)$ such that 
$\phi^*(w)=\phi(w)$ and $\phi^*(w')=\phi'(w')$ for all $w\in W_i,w'\in W'_i, i\in
\mathfrak{r}$. Here if $W_i\cap W'_i\not=\emptyset$ for some $i$
 then we consider a copy of $W'_i$ so that the two domains are disjoint.
\end{df}
\begin{df}[Regularization]\label{070725}
Let $m\ge 0$ and $\varphi
\in \Phi(m).$ 
Let ${\bf G}$ be a $k$-bounded graph on ${\bf \Omega}.$
For an integer $1\le s< k$,
the {\bf $s$-regularization} ${\bf G}/^s\varphi$ is
the $k$-bounded graph on ${\bf \Omega}
$ obtained from ${\bf G}$ by redefining the color of 
each edge ${\bf e}\in {\bf \Omega}_I$ with $I\in {\mathfrak{r}\choose [s]}$
by the 
$\brkt{\sum_{j=0}^{s+1-|I|}\sum_{J\in {\mathfrak{r}\setminus I\choose j}}
m^{j}}$-dimensional 
vector
\begin{eqnarray}
\brkt{{\bf G}/^s\varphi}({\bf e}):=({\bf G}({\bf e}\dot{\cup}{\bf f})
|
J\in {\mathfrak{r}\setminus I\choose [0,s+1-|I|]},
{\bf f}\in {\bf \Omega}_J \, \mbox{\rm with } \, 
{\bf f}\subset \varphi(\mathbb{D})
).\label{081230}
\end{eqnarray}
In the above, when $J=\emptyset,$ we assume ${\bf f}=\emptyset.$
(The sets of colors are naturally extended
while any edge containing at least $s+1$ vertices (i.e. edge of size at least $s+1$)
does not change its (face-)color.)
\par
When $s=k-1$, we simply write ${\bf G}/\varphi:={\bf G}/^{k-1}\varphi.$
\par
For $\vec{\varphi}=(\varphi_i)_{i\in [k-1]}\in\Phi(m_1,\cdots,m_{k-1}),$
we define the {\bf regularization of ${\bf G}$ by $\vec{\varphi}$} by
\begin{eqnarray*}
{\bf G}/\vec{\varphi}:=
(({\bf G}/^{k-1}\varphi_{k-1})/^{k-2}\varphi_{k-2})\cdots/^1\varphi_1.
\end{eqnarray*}
\end{df}
When making ${\bf G}/\vec{\varphi}$ from ${\bf G},$ a size-$s$ edge with $1\le s\le k$ 
 changes its face-color $k-s$ times at the operations 
$/^{k-1}\varphi_{k-1},\cdots,/^s\varphi_s,$ depending on 
$(m_{k-1}+\cdots+m_s)r$ random 
vertices in ${\bf \Omega}.$
 It does not change at the operations 
$/^{s-1}\varphi_{s-1},\cdots,/^1\varphi_1.$
In particular, any size-$k$ (full-size) edge never changes its face-color.
\begin{df}[Regularity]\label{070724a}
Let ${\bf G}$ be a $k$-bounded graph on $ {\bf \Omega}$. 
For
$\vec{\mathfrak{c}}=(\mathfrak{c}_J)_{J\subset I}\in {\rm
TC}_I({\bf G}), I\in {\mathfrak{r}\choose [k]}$, we define {\bf
relative density} by 
\begin{eqnarray*}
{\bf d}_{\bf G}(\vec{\mathfrak{c}}):=
\Prob_{{\bf e}\in {\bf \Omega}_I
}[
{\bf G}({\bf e})=\mathfrak{c}_I
|
{\bf G}(\partial{\bf e})=
(\mathfrak{c}_J)_{J\subsetneq I}
].
\end{eqnarray*}
For a positive
integer $h$ and $\epsilon\ge 0$, we call 
 ${\bf G}$ to be 
 {\bf $(\epsilon,k,h)$-regular}
if and only if there exists a function ${\delta}: {\rm TC}({\bf
G})\to [0,\infty)$ such that
\begin{eqnarray}
\hspace{-5mm}
{\rm (i)}&
 \Prob_{\phi\in\Phi(h)}
[{\bf G}(\phi(e))=S(e)\,, \forall e\in {\mathbb V}(S)]=
\displaystyle\prod_{e\in {\mathbb V}(S)} \brkt{ {\bf d}_{\bf
G}(S\ang{e}) \dot{\pm} \delta(S\ang{e}) },& \forall S\in {\mathcal
S}_{k,h,{\bf G}}, \label{lj19}
\\
\hspace{-5mm} {\rm (ii)}& \Ex_{{\bf e}\in {\bf
\Omega}_I}[\delta({\bf G}\ang{\bf e})]\le \epsilon /|{\rm
C}_I({\bf G})|, & \forall I\in {\mathfrak{r}\choose [k]},
\label{lj26}
\end{eqnarray}
where $a\dot{\pm}b$ means a suitable integer $c$ satisfying $\max\{0,a-b\}\le c\le \min\{1,a+b\}$.
Denote by ${\bf reg}_{k,h}({\bf G})$ the minimum value
of $\epsilon$ such that ${\bf G}$ is $(\epsilon,k,h)$-regular.
\end{df}
The minimum value of $\epsilon$ 
always exists because inequality (\ref{lj26}) includes equality.
Note that if $\delta(\cdot)\equiv 0$ satisfies the above (\ref{lj19})
then the edges of 
${\bf G}$ are colored uniformly at random.
\begin{remark}
Condition (i) measures how far from random
the graph {\bf G} is with respect to containing the expected
number of copies of the (colored) subgraphs $S \in
\mathcal{S}_{h,{\bf G}}$. The smaller $\delta$ is, the closer {\bf
G} is to being random. When $\delta\equiv 0$, then {\bf G} behaves
exactly like a random graph. On the other hand, if we take $\delta
\equiv 1$ then (i) is automatically satisfied. Condition (ii) places an
upper bound on the size of $\delta$.
Our proof will yield the main theorem even if we replace
the right-hand side of (ii) by $g_I(|{\rm C}_I({\bf G})|)$
for any fixed functions $g_I>0$, for example,
$g_I(x)=x^{-1/\epsilon}.$
\end{remark}
\begin{remark}
In $\Prob_{{\bf e}\in {\bf \Omega}_I}
[\cdots
]$ and 
$\Ex_{{\bf e}\in {\bf \Omega}_I}
[\cdots
]$, ${\bf e}$ is a random variable, equivalently a sequence of 
$|I|$ random vertices.
The relative density 
${\bf d}_{\bf G}(\vec{\mathfrak{c}})$
 is undefined when 
$\Prob_{{\bf e}\in {\bf \Omega}_I}
[{\bf G}(\partial {\bf e})=
(\mathfrak{c}_J)_{J\subsetneq I}]=0.$
But this will not cause any trouble later, in particular at 
(\ref{lj19}), since such a relative density will be always multiplied by zero.
Here we define ${\bf d}_{\bf G}(\vec{\mathfrak{c}})$
 to be one, if 
 $\Prob_{{\bf e}\in {\bf \Omega}_I}
[{\bf G}(\partial {\bf e})=
(\mathfrak{c}_J)_{J\subsetneq I}]=0.$
\end{remark}
\par
Our main theorem is as follows.
\begin{thm}[Main Theorem]\label{070908b}
For any $r\ge k,\, h,\, \vec{b}=(b_i)_{i\in [k]},$ and $\epsilon>0,$
there exist (increasing) functions $
m^{(i)}:\naturalset^{k-i}\to\naturalset $
 and
$\tilde{n}^{(i)} :\naturalset^{k-1-i}\to \naturalset,\, i\in
[k-1]$ satisfying the following:
\par
If
 ${\bf G}$  is
a $\vec{b}$-colored ($k$-bounded $r$-partite hyper)graph
on $
{\bf \Omega}$ then we have 
\begin{eqnarray*}
\Ex_{\vec{n}=(n^{(1)},\cdots,n^{(k-1)})}
\Ex_{\vec{\varphi}=(\varphi_i)_{i\in [k-1]}}
[{\bf reg}_{k,h}({\bf G}/\vec{\varphi})]
\le \epsilon.
\end{eqnarray*}
In the above probabilistic process, 
each integer $n^{(i)}$ (from $i=k-1$ to $i=1$)
 is picked uniformly at random 
from $[0,\tilde{n}^{(i)}(n^{(i+1)},\cdots,
n^{(k-1)})-1]$.
Each $\varphi_i\in \Phi(m^{(i)}(n^{(i)},\cdots,n^{(k-1)}))
$ is random.
\end{thm}
In the above, $\tilde{n}^{(k-1)}$ is read to be a constant integer.
When $k=1$, the theorem is read to be true trivially where 
we do not take $\vec{n}$ and put ${\bf G}/\vec{\varphi}={\bf G}$ 
while any $1$-bounded ${\bf G}$ is $(0,1,h)$-regular.
Thus ${\bf reg}_{1,h}({\bf G})=0.$
\par
Note that $m^{(i)},\tilde{n}^{(i)}$ depend only on
$r,k,h,\vec{b},\epsilon$ and are independent of everything else
including ${\bf \Omega}$. The following immediate consequence is convenient for applications.
\\
\begin{cor}[Regularity Lemma (including so-called Counting Lemma)]
\label{lj21} For any $r\ge k,\,h,\,\vec{b}=(b_i)_{i\in
[k]},\,\epsilon>0,$ there exist integers $
\widetilde{m}_1,\cdots,\widetilde{m}_{k-1} $ such that if
 ${\bf G}$  is
a $\vec{b}$-colored ($k$-bounded $r$-partite hyper)graph on $
{\bf \Omega}$ then  for some
integers $m_1,\cdots,m_{k-1}$ with
$m_i\le \widetilde{m}_i, i\in [k-1],$
\begin{eqnarray}
\Ex_{\vec{\varphi}\in\Phi(m_1,\cdots,m_{k-1})}
[{\bf reg}_{k,h}({\bf G}/\vec{\varphi})]
\le \epsilon.\label{081215}
\end{eqnarray}
In particular, when (\ref{081215}) holds, if we pick a map 
$\vec{\varphi}\in\Phi(m_1,\cdots,m_{k-1})$ randomly then with probability at least 
$1-\sqrt{\epsilon}$, we have ${\bf reg}_{k,h}({\bf G}/\vec{\varphi})\le \sqrt{\epsilon},$ 
thus 
${\bf G}/\vec{\varphi}$ is $(\sqrt{\epsilon},k,h)$-regular.
\end{cor}
\begin{example}
If $r=k=h=2$ and $(b_1,b_2)=(1,2)$ then the corollary becomes 
 one of the usual Graph Regularity Lemmas, when 
 ${\bf G}$  has black and white edges and 
$S$ is an ordinary bipartite graph on 
$\{u_1,v_1\}\dot{\cup}\{u_2,v_2\}$ such that 
$u_1$ and $v_1$ have the same color, say red$_1$, 
that $u_2$ and $v_2$ have the same color, say red$_2$, 
and that the four 
edges $u_1u_2,u_1v_2,v_1u_2,v_1v_2$ have the same color, say black.
(The color red$_i$ may be considered as 
a sequence of black and white colors.)
\end{example}
\hspace{5mm} Our proof will yield the theorem even if we replace
the right-hand side of (\ref{lj26}) by $g_I(|{\rm C}_I({\bf G})|)$
for any fixed functions $g_I>0$, for example,
$g_I(x)=x^{-1/\epsilon}.$
If the reader is interested only in applications 
to Szemer\'{e}di's theorem, then 
it suffices to consider only the case of $h=1.$
\section{Proof of the Main Theorem}
\subsection{Two lemmas and their proofs}
\begin{df}[Notation for the lemmas]\label{081229}
Let ${\bf G}$ be an ($r$-partite $(b_i)_{i\in [k]}$-colored) 
$k$-bounded graph on ${\bf \Omega}$.
For two edges ${\bf e}, {\bf e'}\in {\bf \Omega}_I,$ we abbreviate
$
{\bf G}({\bf e})
=
{\bf G}({\bf e'})$ and $
{\bf G}(\partial {\bf e})
=
{\bf G}(\partial {\bf e'})
$
by 
$
{\bf e}\stackrel{\bf G}{\approx}{\bf e'}$ and $
{\bf e}\stackrel{\partial {\bf G}}{\approx}{\bf e'},
$ respectively.
\par
An {\bf $(s,h)$-error function} of ${\bf G}$ is
a function $\delta:\bigcup_{I\in {\mathfrak{r}\choose [s]}}{\rm TC}_I({\bf G})\to
[0,\infty)$
satisfying (\ref{lj19}) for all
$S\in {\mathcal S}_{s,h,{\bf G}}$.
We write 
${\bf d}_{\bf G}^{(\delta)}(\vec{\mathfrak{c}})
:={\bf d}_{\bf G}(\vec{\mathfrak{c}})\dot{\pm}\delta(\vec{\mathfrak{c}})$
 and 
${\bf d}_{\bf G}^{(-\delta)}(\vec{\mathfrak{c}})
:={\bf d}_{\bf G}(\vec{\mathfrak{c}})-\delta(\vec{\mathfrak{c}})$
 for
$\vec{\mathfrak{c}}\in {\rm TC}({\bf G}).$
\par 
We abbreviate $\bigcup_{i\in [k-1]}{\mathbb V}_i(S)$ by ${\mathbb
V}_{(k-1)}(S)$.
\par
Denote by $\dbrkt{\cdots}$ the Iverson bracket, i.e., it equals
$1$ if the statement in the bracket holds, and $0$ otherwise.
\end{df}
\begin{lemma}[Correlation bounds counting error]\label{lj23c}
For a $k$-bounded graph ${\bf G}$ and $S\in {\mathcal S}_{k,h,{\bf G}}$,
we have that
\begin{eqnarray*}
&& \left| \Prob_{\phi\in\Phi(h)}\left[\left. {\bf
G}(\phi(e))=S(e), \,\forall e\in {\mathbb V}_k(S) \right| {\bf
G}(\phi(e))=S(e),\,\forall e\in {\mathbb V}_{(k-1)}(S)
\right]\nonumber - \prod_{e\in {\mathbb V}_k(S)} {\bf d}_{\bf
G}(S\ang{e}) \right|
\\
&\le& |{\mathbb V}_k(S)| \max_{\emptyset\not=D\subset {\mathbb
V}_k(S)} \left| \Ex_{\phi\in\Phi(h)}\left[\left. \prod_{e\in D}
\brkt{ \dbrkt{{\bf G}(\phi(e))=S(e)} - {\bf d}_{\bf G}(S\ang{e}) }
\right| {\bf G}(\phi(e))=S(e),\,\forall e\in {\mathbb
V}_{(k-1)}(S) \right] \right| . \label{lj15b}
\end{eqnarray*}
\end{lemma}
\Proof We prove it by induction on $|{\mathbb V}_k(S)|.$ If
$|{\mathbb V}_k(S)|=0$ or $1$ then it is trivial, since 
in this case, the left-hand side of the inequality is $0.$
So let us assume that $|{\mathbb V}_k(S)|\ge 2$ and that 
the result holds for all smaller values of $|\mathbb{V}_k(S)|.$
Let $d_e:=
{\bf d}_{\bf G}(S\ang{e})$ and let $\eta$ be the maximum part of
the desired right-hand side. 
Then for $D:=\mathbb{V}_k(S)$ we have 
\begin{eqnarray*}
[-\eta,\eta]&\ni&
 \Ex_{\phi\in\Phi(h)=\Phi(V(S))}\left[\left. \prod_{e\in {\mathbb
V}_k(S)}
\brkt{
\dbrkt{{\bf G}(\phi(e))=S(e)} - {\bf d}_{\bf G}(S\ang{e})
}
\right| {\bf G}(\phi(e))=S(e),\,\forall e\in {\mathbb V}_{(k-1)}(S)
\right]\nonumber\\
&=&  \Ex_{\phi\in\Phi(h )}\left[\left. \prod_{e\in {\mathbb
V}_k(S)} \dbrkt{{\bf G}(\phi(e))=S(e) } \right| {\bf
G}(\phi(e))=S(e),\,\forall e \in {\mathbb V}_{(k-1)}(S) \right]\nonumber
\\
&&+ \sum_{\emptyset\not=D\subset {\mathbb V}_k(S)}
\brkt{\prod_{e\in D} (-d_e) } \Ex_{\phi\in\Phi(h) }\left[ \left.
\prod_{e\in {\mathbb V}_k(S)\setminus D} \dbrkt{{\bf
G}(\phi(e))=S(e) } \right| {\bf G}(\phi(e))=S(e),\forall e \in
{\mathbb V}_{(k-1)}(S) \right], \label{081021}
\end{eqnarray*}
expanding the product and using the
linearity of expectation and the definition of $d_e$.
Now we will focus on second term above. Since the value of
$\dbrkt{{\bf G}(\phi(e))=S(e) }$ is $0$ or $1$, we can replace
$\Ex$ by $\Prob$, and consequently, apply the induction hypothesis
(since $D$ is nonempty). 
Consider a complex $S^-$ with $\mathbb{V}_k
(S^-)=\mathbb{V}_k(S)\setminus D$ by 
invisualizing the edges in $D$ of $S$. 
\\[.1in]
Using the inductive hypothesis for complex $S^-$ 
 in the place of $S$, we rewrite the second term and obtain 
\begin{eqnarray*}
&& \Ex_{\phi\in\Phi(h )}\left[\left. \prod_{e\in {\mathbb V}_k(S)}
\dbrkt{{\bf G}(\phi(e))=S(e) } \right| {\bf
G}(\phi(e))=S(e),\,\forall e\in {\mathbb V}_{(k-1)}(S) \right]
\nonumber
\\
&
\stackrel{\rm I.H.}
{=}
& - \sum_{\emptyset\not=D\subset {\mathbb
V}_k(S)} \brkt{\prod_{e\in D} (-d_e) } \brkt{ 
\brkt{\prod_{e\in {\mathbb
V}_k(S)\setminus D} d_e}\pm
|{\mathbb V}_k(S^-)|
 \eta
 } \quad 
\dot{\pm}\quad \eta
\nonumber
\\
&=& 
- \brkt{\brkt{ \prod_{e\in {\mathbb V}_k(S)} d_e}\pm
|{\mathbb V}_k(S^-)| \eta
} \sum_{\emptyset\not=D\subset {\mathbb V}_k(S)} \brkt{\prod_{e\in
D} (-1) } \quad \dot{\pm}\, \eta \quad (
\because |d_e|\le 1)
\nonumber
\\
&=& 
- \brkt{ 
\brkt{\prod_{e\in {\mathbb V}_k(S)} d_e}\pm
(|{\mathbb V}_k(S)|-1) \eta
}
\brkt{ (1-1)^{|{\mathbb V}_k(S)|} -1 } 
\quad \dot{\pm}\,\eta 
 \quad (\because 
|\mathbb{V}_k(S)|>|\mathbb{V}_k(S^-)|
)
\nonumber
\\
&=&
\brkt{
\prod_{e\in {\mathbb V}_k(S)}
d_e}\pm|{\mathbb V}_k(S)|\eta.
\end{eqnarray*}
\lemmaqed
We will use the following form of the Cauchy-Schwarz.
\begin{fact}[Cauchy-Schwarz inequality]\label{080827a}
For a random variable $X$ on a probability space ${\Omega}$ if 
an equivalent relation $\approx$ on $\Omega$ is a refinement of another equivalent relation $\sim$ 
on $\Omega$ then
\begin{eqnarray}
\Ex_{{\omega_0}\in {\Omega}}\brkt{\Ex_{\omega\in \Omega}
[X({\omega})|{\omega}\approx {\omega_0}]}^2
\ge 
\Ex_{\omega_0\in\Omega}\brkt{\Ex_{\omega\in\Omega}[X({\omega})| {\omega}\sim {\omega_0}]
}^2.\label{080827}
\end{eqnarray}
\end{fact}
\Proof 
By the Cauchy-Schwarz  
(i.e. $\Ex[X^2]\Ex[Y^2]\ge (\Ex[XY])^2$), 
we have 
$
\Ex_{\omega_0}\brkt{\Ex_{\omega}[X({\omega})|{\omega}\approx {\omega_0}]}^2
=
\Ex_{\omega_0}\left[
\Ex_{\omega'}\left[\left.
\brkt{\Ex_{\omega}[X({\omega})|{\omega}\approx {\omega'}]
}^2\right|{\omega'}\sim {\omega_0}
\right]\right]
=
\Ex_{\omega_0}\left[
\Ex_{\omega'}[1^2|{\omega'}\sim {\omega_0}]
\cdot
\Ex_{\omega'}\left[\left.
\brkt{\Ex_{\omega}[X({\omega})|{\omega}\approx {\omega'}]
}^2\right|{\omega'}\sim {\omega_0}
\right]\right]
\stackrel{\rm CS}{\ge}
\Ex_{\omega_0}\left(
\Ex_{\omega'}[
{1\cdot\Ex_{\omega}[X({\omega})|{\omega}\approx {\omega'}]
}|{\omega'}\sim {\omega_0}
]\right)^2= 
\Ex_{\omega_0}\brkt{\Ex_{\omega}[X({\omega})| {\omega}\sim {\omega_0}]
}^2.
$
\factqed
With this fact and Definition \ref{081229}, we next tackle
\begin{lemma}[Mean square bounds correlation]
\label{alj08}
Let $k,h,m$ be positive integers and
${\bf G}$ a $k$-bounded graph on vertex sets ${\bf \Omega}$.
Let $S\in {\mathcal S}_{k,h,{\bf G}}$ and let
$F_e:{\rm C}_I({\bf G})
\to [-1,1]$ be a function for each
$I\in {\mathfrak{r}\choose k}$ and
for each $e\in {\mathbb V}_I(S).$
If $\delta$ is a
$(k-1,2h)$-error function of ${\bf G}$ then
for any $I\in {\mathfrak{r} \choose k}$ and
$e_0\in {\mathbb V}_I(S),$ we have that
\\
\begin{eqnarray}
&&
\brkt{
\Ex_{\phi\in\Phi(h)}\left[
\prod_{e\in {\mathbb V}_k(S)}
F_e\brkt{\phi(e)}
\prod_{e\in {\mathbb V}_{(k-1)}(S)}
\llbracket
{\bf G}(\phi(e))=S(e)
\rrbracket
\right]}^2
\nonumber\\
&\le&
\Ex_{\varphi\in \Phi(mh)}
\Ex_{{\bf e^*}\in {\bf \Omega}_I
}
[
\brkt{
\Ex_{{\bf e}\in {\bf \Omega}_I
}[F_{e_0}\brkt{\bf e}
\llbracket
{\bf G}(\partial {\bf e})=S(\partial e_0)\rrbracket
|
{\bf e}\stackrel{
\partial {\bf G}/\varphi}{\approx}{\bf e^*}]
}^2
]\nonumber
\\
&&\cdot
\brkt{
\prod_{e\in {\mathbb V}_{(k-1)}(S)}
{\bf d}_{\bf G}^{(\delta)}(S\ang{e})
}
\brkt{
\brkt{
\prod_{e\in {\mathbb V}_{(k-1)}(S),e\not\subset e_0}
{\bf d}_{\bf G}^{(\delta)}(S\ang{e})
}
+
{
1 \over m}
}
\label{070907d}
\end{eqnarray}
where $\phi, \varphi$ are random and where we abbreviate 
 $F_e({\bf G}({\bf e}))$ 
by $F_e({\bf e})$ (thus, $F_e(\phi(e))=F_e({\bf G}(\phi(e)))$).
\par
In particular, if we suppose 
\begin{eqnarray}
\min_{
J\in {\mathfrak{r}\choose [k-1]}}
\min_{e\in {\mathbb V}_J(S)}
\brkt{
 {1\over 2}
{\bf d}_{\bf G}(S\ang{e}) 
-
\delta(S\ang{e})
}
>0\, 
\mbox{ and } \, 
{1\over m}
\le
\prod_{e\in {\mathbb V}_{(k-1)}(S), e\not\subset e_0
}
{\bf d}_{\bf G}^{(-\delta)}(S\ang{e}) \label{070907b}
\end{eqnarray}
(i.e. $\delta$ is small and $m$ is large) then
\begin{eqnarray}
&&
\brkt{
\Ex_{\phi\in\Phi(h)}\left[
\left.
\prod_{e\in {\mathbb V}_k(S)}
F_e\brkt{\phi(e)}
\right|
{\bf G}(\phi(e))=S(e) \, \forall 
e\in {\mathbb V}_{(k-1)}(S)
\right]
}^2\nonumber
\\
&\le& 
2\cdot 3^{2|\mathbb{V}_{(k-1)}(S)|}
\Ex_{\varphi\in \Phi(mh)}
\Ex_{{\bf e^*}\in {\bf \Omega}_I
}
[
\brkt{
\Ex_{{\bf e}\in {\bf \Omega}_I
}[F_{e_0}\brkt{\bf e}
|
{\bf e}\stackrel{
\partial {\bf G}/\varphi}{\approx}{\bf e^*}]
}^2
|{\bf G}(\partial {\bf e^*})=S(\partial e_0)
].\label{070907c}
\end{eqnarray}
\end{lemma}
\Proof [Tools: Cauchy-Schwarz, Fact \ref{080827a}] Fix $I_0\in
{\mathfrak{r}\choose k}$ and $e_0\in {\mathbb V}_{I_0}(S)$. For
$\phi\in\Phi(V(S)\setminus e_0)$ and for ${\bf e_0}\in {\bf
\Omega}_{I_0}$, we define the (extended) function $\phi^{({\bf
e_0})}\in\Phi(V(S))$ such that:\\
 (i) each $v\in e_0$ is mapped to
the corresponding ${\bf v}\in {\bf e_0}$ with the index of $v$,
(thus, if $v\in e_0$ has an index $i\in\mathfrak{r}$
then $
{\bf e_0}\cap {\bf \Omega}_i=\{
\phi^{({\bf e_0})}(v)\}$
)
and
 that,\\
  (ii) each $v\in V(S)\setminus e_0$ is mapped to
$\phi(v)$.\\ 
(That is, when we have a map $\phi$ defined for $r-k$ vertices, 
we extend it by assigning the remaining $k$ vertices in $e_0$ to 
the $k$ vertices in ${\bf e_0}$ so that it will be a 
partitionwise map from $V(S)$ to ${\bf \Omega}$.)
For $\vec{\varphi}=(\varphi_i)_{i\in [m]}$ with
$\varphi_i\in\Phi(V(S)\setminus e_0),$ we define an equivalence
relation $\stackrel{\vec{\varphi}}{\sim}$ on ${\bf \Omega}_{I_0}$
by the condition that
\begin{eqnarray}
{\bf e}\stackrel{\vec{\varphi}}{\sim}{\bf e'} \mbox{ if and only
if } \varphi_i^{({\bf e})}(e) \stackrel{\bf G}{\approx}
\varphi_i^{({\bf e'})}(e), \quad  \forall e\in
\mathbb{V}(S)\setminus\{e_0\}, \forall i\in [m].\label{lj23}
\end{eqnarray}
(Note that $V(S)\setminus e_0$ is a vetex set while 
$\mathbb{V}(S)\setminus \{e_0\}$ is an edge set. 
Since the right-hand side of (\ref{lj23}) holds trivially for $e$ with 
 $e\cap e_0=\emptyset$, it is enough to check only for 
$e$ with $1\le |e\cap e_0|\le k-1.$
)\\
Let $S^{(1)},\cdots,S^{(m)}$ and $e_0^{(1)},\cdots,e_0^{(m)}$ be
copies of $S$ and of $e_0$. For $\vec{\varphi}=(\varphi_i)_{i\in
[m]}$ with $\varphi_i\in\Phi (V(S^{(i)})\setminus e_0^{(i)})$ if
$\varphi^*\in
\Phi(mh)=\Phi(V(S^{(1)})\dot{\cup}\cdots\dot{\cup}V(S^{(m)}))$ is
an extended function of $\varphi_i$'s, i.e.,
$\varphi^*(v)=\varphi_i(v)$ for any $v\in V(S^{(i)})\setminus
e_0^{(i)}, i\in [m]$ then, because of (\ref{081230}) and 
(\ref{lj23}), 
 it is easily seen that
\begin{eqnarray}
{\bf e}\stackrel{\partial {\bf G}/\varphi^*}{\approx}{\bf e'}
\mbox{ implies }
{\bf e}\stackrel{\vec{\varphi}}{\sim}{\bf e'}
\label{lj23a}
\end{eqnarray}
where ${\bf G}/\varphi^*={\bf G}/^{k-1}\varphi^*
$ is the $(k-1)$-regularization.
\\
(To see this, observe that 
$
{\bf e}\stackrel{\partial {\bf G}/\varphi^*}{\approx}{\bf e'}
$ means that ${\bf G}/\varphi^*({\bf e}|_J)
=
{\bf G}/\varphi^*({\bf e'}|_J)
$ for all $J\subsetneq I_0.$ By (\ref{081230}), if 
$J'\in {\mathfrak{r}\setminus J\choose [0,k-|J|]}, 
{\bf f}\in {\bf \Omega}_{J'}$ and ${\bf f}\subset \varphi^*(\mathbb{D})
$ then ${\bf e}|_J\dot{\cup}{\bf f}
\stackrel{\bf G}{\approx}
{\bf e'}|_J\dot{\cup}{\bf f}
$. Since 
$\left|{\bf e}|_J\dot{\cup}{\bf f}\right|\le k$, for all 
$e\in \mathbb{V}(S)\setminus\{e_0\}$
we have $\varphi_i^{({\bf e})}(e)
=
{\varphi^{*}}^{({\bf e})}_i(e)
\stackrel{\bf G}{\approx}
{\varphi^{*}}^{({\bf e'})}_i(e)
=
\varphi_i^{({\bf e'})}(e),
$ where ${\varphi^*_i}^{({\bf e})}$ and 
 ${\varphi^*_i}^{({\bf e'})}$ are naturally defined by 
restricting the domain of $\varphi_i^*$ 
from $V(S^{(i)})$ 
to $V(S^{(i)})
\setminus e_0.$
By (\ref{lj23}), ${\bf e}\stackrel{\vec{\varphi}}{\sim}{\bf e'}$.)
\par
Let $F^*_{e_0}({\bf e}):=
F_{e_0}({\bf e})\dbrkt{{\bf G}(\partial {\bf e})=S(\partial e_0)}
$ and let 
\begin{eqnarray*}
F^*(\phi):=
\prod_{e\in {\mathbb V}_k(S)\setminus
\{e_0\}}
F_e\brkt{\phi(e)}\prod_{e\in {\mathbb V}_{(k-1)}(S)}
\dbrkt{{\bf G}(\phi(e))=S(e)}.
\end{eqnarray*}
Note that ${\bf G}(\partial\phi(e_0))=S(\partial e_0)
$ holds if-and-only-if 
${\bf G}(\phi(e))=S(e)
$ for all $e\subsetneq e_0.$
Also $\dbrkt{P}\in \{0,1\}$ implies 
$\dbrkt{P}^2=\dbrkt{P}$ for any statement $P.$
With the two facts, the left-hand side of (\ref{070907d})
equals 
\begin{eqnarray}
&&
\brkt{
\Ex_{\phi\in\Phi(h)}\left[
\prod_{e\in {\mathbb V}_k(S)}
F_e\brkt{\phi(e)}
\prod_{e\in {\mathbb V}_{(k-1)}(S)}
\llbracket
{\bf G}(\phi(e))=S(e)
\rrbracket
\right]}^2
\nonumber\\
&=&
\brkt{
\Ex_{\phi}\left[
F_{e_0}(\phi(e_0))
\prod_{e\in {\mathbb V}_k(S)\setminus\{e_0\}}
F_e\brkt{\phi(e)}\cdot
\dbrkt{{\bf G}(\partial \phi(e_0))=S(\partial e_0)}^2
\prod_{e\in {\mathbb V}_{(k-1)}(S):e\not\subset e_0}
\llbracket
{\bf G}(\phi(e))=S(e)
\rrbracket
\right]}^2
\nonumber\\
&=&
\brkt{
\Ex_{\phi\in\Phi(h)}\left[
F^*_{e_0}\brkt{\phi(e_0)}F^*(\phi)
\right]}^2 
\quad\quad\quad (\mbox{by definitions of }F^*_{e_0} \mbox{ and }F^*)
\nonumber
\\
&=&
\brkt{
\Ex_{{\bf e_0}\in {\bf \Omega}_{I_0},
\phi\in\Phi(V(S)\setminus e_0)}\left[
F^*_{e_0}\brkt{\bf e_0}F^*(\phi^{({\bf e_0})})
\right]
}^2
\nonumber\\
&& 
\quad\quad\quad 
(\mbox{since $\phi\in\Phi(h)=\Phi(V(S))$ 
consists of $rh=k+(rh-k)$ random vertices in ${\bf \Omega}$})
\nonumber\\
&=&
\brkt{\Ex_{\vec{\varphi}=(\varphi_i)_{i\in [m]}\in (\Phi(V(S)\setminus e_0))^m}
\Ex_{{\bf e_0}\in {\bf \Omega}_{I_0}}
\left[
F^*_{e_0}\brkt{\bf e_0}
\Ex_{i\in [m]}[F^*(\varphi_i^{({\bf e_0})})]
\right]
}^2
\nonumber\\
&=&
\brkt{\Ex_{\vec{\varphi}=(\varphi_i)_{i\in [m]}\in (\Phi(V(S)\setminus e_0))^m}
\Ex_{{\bf e_0}\in {\bf \Omega}_{I_0}}
\left[
\Ex_{{\bf e}\in {\bf \Omega}_{I_0}}
\left[
F^*_{e_0}\brkt{\bf e}
\Ex_{i\in [m]}[F^*(\varphi_i^{({\bf e})})]
\left|{\bf e}\stackrel{\vec{\varphi}}{\sim}{\bf e_0}\right.
\right]
\right]
}^2
\nonumber\\
&
\stackrel{(\ref{lj23})}{=}&
\brkt{\Ex_{\vec{\varphi}=(\varphi_i)_{i\in [m]}\in (\Phi(V(S)\setminus e_0))^m}
\Ex_{{\bf e_0}\in {\bf \Omega}_{I_0}}
\left[
\Ex_{{\bf e}\in {\bf \Omega}_{I_0}}
[
F^*_{e_0}\brkt{\bf e}
|{\bf e}\stackrel{\vec{\varphi}}{\sim}{\bf e_0}
]
\Ex_{i\in [m]}[F^*(\varphi_i^{({\bf e_0})})]
\right]
}^2
\nonumber\\
&&
\mbox{\hspace{1.5in} (since $F^*(\varphi^{({\bf
e})}_i)=F^*(\varphi^{({\bf e_0})}_i)$ by (\ref{lj23}) when ${\bf
e}\stackrel{\vec{\varphi}}{\sim}{\bf e_0}$)}
\nonumber\\
&\le^{\rm C.S.}&
{
\Ex_{\vec{\varphi}
}
\Ex_{{\bf e_0}}
\left[
\brkt{
\Ex_{{\bf e}\in {\bf \Omega}_{I_0}}
[
F^*_{e_0}\brkt{\bf e}
| {\bf e}\stackrel{\vec{\varphi}}{\sim}{\bf e_0}
]
}^2
\right]
\cdot
\Ex_{\vec{\varphi}=(\varphi_i)_i
}
\Ex_{{\bf e_0}}
\left[
\brkt{
\Ex_{i\in [m]}[
F^*(\varphi_i^{({\bf e_0})})
]
}^2
\right]
}
\nonumber\\
&\stackrel{(\ref{lj23a}),(\ref{080827})}{\le}
&
\Ex_{\varphi^*\in\Phi(mh)
}
\Ex_{{\bf e_0}}
\left[
\brkt{
\Ex_{{\bf e}\in {\bf \Omega}_{I_0}}
[
F^*_{e_0}({\bf e})
| {\bf e}\stackrel{\partial {\bf G}/\varphi^*}{\approx}{\bf e_0}
]
}^2
\right]
\cdot
\Ex_{\vec{\varphi}=(\varphi_i)_i
}\Ex_{{\bf e_0}\in {\bf \Omega}_{I_0}}
\left[
\Ex_{i,j\in [m]}[
F^*(\varphi_i^{({\bf e_0})})
F^*(\varphi_j^{({\bf e_0})})
]
\right]
.
\nonumber
\end{eqnarray}
The first term 
of the last line appears in the first term of our desired upperbound. 
We now forcus on the second term. Since $|F^*(\cdot)|\le 1,$ it 
equals
\begin{eqnarray}
&&
\Ex_{{\bf e_0}\in {\bf \Omega}_{I_0}}
\left[
{m-1\over m}
\Ex_{i\not=j\in [m]}
\Ex_{\vec{\varphi}=(\varphi_i)_i}
[
F^*(\varphi_i^{({\bf e_0})})
F^*(\varphi_j^{({\bf e_0})})
]
+
{1\over m}
\Ex_{i\in [m]}
\Ex_{\vec{\varphi}=(\varphi_i)_i}[
\brkt{
F^*(\varphi_i^{({\bf e_0})})
}^2
]
\right]\nonumber
\\
&\le &
\Ex_{{\bf e_0}\in {\bf \Omega}_{I_0}}
\Ex_{\varphi_1,\varphi_2\in \Phi(V(S)\setminus e_0)}
[
F^*(\varphi_1^{({\bf e_0})})
F^*(\varphi_2^{({\bf e_0})})
]
+
{1\over m}\Ex_{{\bf e_0}\in {\bf \Omega}_{I_0}}
\Ex_{\varphi_1\in \Phi(V(S)\setminus e_0)}
[
\left|F^*(\varphi_1^{({\bf e_0})})\right|
]
\nonumber\\
&\le &
\Ex_{{\bf e_0}\in {\bf \Omega}_{I_0}}
\Ex_{\varphi_1,\varphi_2\in \Phi(V(S)\setminus e_0)}
\left[
\prod_{e\in {\mathbb V}_{(k-1)}(S)}
\dbrkt{{\bf G}(\varphi_1^{({\bf e_0})}(e))=S(e)}
\dbrkt{{\bf G}(\varphi_2^{({\bf e_0})}(e))=S(e)}
\right]
\nonumber\\
&&
+
{1\over m} \Ex_{\varphi\in \Phi(V(S))} \left[ \prod_{e\in
{\mathbb V}_{(k-1)}(S)} \dbrkt{{\bf G}(\varphi(e))=S(e)} \right].
(\mbox{by the definition of }F^* \mbox{ since }|F_e|\le 1)
\label{081231}
\end{eqnarray}
Looking at the second term first, this can be written as
\begin{eqnarray}
{1\over m} \Prob_{\varphi\in \Phi(V(S))} 
\left[ 
{{\bf G}(\varphi(e))=S(e)} \forall e\in{\mathbb V}_{(k-1)}(S)
\right]
=
{1\over m} \prod_{e\in {\mathbb
V}_{(k-1)}(S) } {\bf d}^{(\delta)}_{\bf G}(S\ang{e}),
\label{081231a}
\end{eqnarray}
applying the assumption that $\delta$ is $(k-1,2h)$-error function of ${\bf G}$ 
to an $S^-\in {\mathcal S}_{k-1,h,{\bf G}}$ with $\mathbb{V}(S^-):=\mathbb{V}_{(k-1)}(S)$.
\par
We will interpret the first term by applying the same assumption on $\delta$ to 
another complex $S''.$
Here $S''\in \mathcal{S}_{k-1,2h,{\bf G}}$ is a simplicial-complex 
obtained from two copies of $S^-$, say ${S^-}^{(1)}$ and ${S^-}^{(2)},$
by identifying any pair of vertices 
$v^{(1)}\in e_0^{(1)}$ 
and $v^{(2)}\in e_0^{(2)}$ in which 
$e_0^{(1)}$ and $e_0^{(2)}$ are the edges in the copies of 
$S^-$ corresponding to $e_0$.
(Any edge $e$ containing two vertices $v^{(1)}\in V({S^-}^{(1)})\setminus e_0^{(1)}$ and 
$v^{(2)}\in V({S^-}^{(2)})\setminus e_0^{(2)}$ is invisible in $S''$.)
Applying the assumption on $\delta$ to this $S'',$ the first term can be rewrriten as
\begin{eqnarray*}
&& 
\Ex_{{\bf e_0}}
\Ex_{\varphi_1,\varphi_2}
\left[
\prod_{e\in {\mathbb V}_{(k-1)}(S),e\not\subset e_0}
\dbrkt{{\bf G}(\varphi_1^{({\bf e_0})}(e))=S(e)}
\dbrkt{{\bf G}(\varphi_2^{({\bf e_0})}(e))=S(e)}
\prod_{e\subsetneq e_0}
\dbrkt{{\bf G}(\varphi_1^{({\bf e_0})}(e))=S(e)}
\right] 
\nonumber\\ 
&=&
\Ex_{\phi\in \Phi(V(S''))}[
\prod_{e\in {\mathbb V}(S'')} 
\dbrkt{{\bf G}(\varphi(e))=S''(e)}
]\\
&=&\prod_{e\in {\mathbb V}(S'')}
{\bf
d}^{(\delta)}_{\bf G}(S''\ang{e})
\quad\quad\quad (\mbox{since $S''\in {\mathcal S}_{k-1,2h}$ and 
$\delta$ is a $(k-1,2h)$-error function of ${\bf G}$})
\\
&=& 
\prod_{e\in {\mathbb V}_{(k-1)}(S), e\not\subset e_0 } ({\bf
d}^{(\delta)}_{\bf G}(S\ang{e}))^2 \prod_{e\subsetneq e_0} 
{\bf
d}^{(\delta)}_{\bf G}(S\ang{e}),
\end{eqnarray*}
completing the proof of (\ref{070907d}) by (\ref{081231}) and (\ref{081231a}).
\par
Next, we show the last sentence of the lemma.
The left-hand side of (\ref{070907c}) is at most 
\begin{eqnarray*}
&&
\brkt{
\Ex_{\phi\in\Phi(h)}\left[
\prod_{e\in {\mathbb V}_k(S)}
F_e\brkt{\phi(e)}
\prod_{e\in {\mathbb V}_{(k-1)}(S)}
\llbracket
{\bf G}(\phi(e))=S(e)
\rrbracket
\right]
/
\Prob_{\phi\in\Phi(h)}\left[
{\bf G}(\phi(e))=S(e) \, \forall 
e\in {\mathbb V}_{(k-1)}(S)
\right]
}^2\nonumber
\\
&=&\brkt{
\Ex_{\phi\in\Phi(h)}\left[
\prod_{e\in {\mathbb V}_k(S)}
F_e\brkt{\phi(e)}
\prod_{e\in {\mathbb V}_{(k-1)}(S)}
\llbracket
{\bf G}(\phi(e))=S(e)
\rrbracket
\right]
/
\brkt{
\prod_{e\in {\mathbb V}_{(k-1)}(S)}
{\bf d}_{\bf G}^{(\delta)}(S\ang{e})
}}^2
\\
&\stackrel{(\ref{070907d})}{\le }&
\Ex_{\varphi\in \Phi(mh)}
\Ex_{{\bf e^*}\in {\bf \Omega}_I
}
[
\brkt{
\Ex_{{\bf e}\in {\bf \Omega}_I
}[F_{e_0}\brkt{\bf e}
\llbracket
{\bf G}(\partial {\bf e})=S(\partial e_0)\rrbracket
|
{\bf e}\stackrel{
\partial {\bf G}/\varphi}{\approx}{\bf e^*}]
}^2
|{\bf G}(\partial {\bf e^*})=S(\partial e_0)
]
\\
&& 
\cdot
\Prob_{\bf e^*}[{\bf G}(\partial {\bf e^*})=S(\partial e_0)]
\brkt{
\prod_{e\in {\mathbb V}_{(k-1)}(S)}
{\bf d}_{\bf G}^{(\delta)}(S\ang{e})
}
\brkt{
\brkt{
\prod_{e\in {\mathbb V}_{(k-1)}(S),e\not\subset e_0}
{\bf d}_{\bf G}^{(\delta)}(S\ang{e})
}
+
{
1 \over m}
}
\\&&
/
\brkt{
\prod_{e\in {\mathbb V}_{(k-1)}(S)}
{\bf d}_{\bf G}^{(-\delta)}(S\ang{e})
}^2 
\\
&&
(
\mbox{since ${\bf e}\stackrel{\partial {\bf G}/\varphi}{\approx}{\bf e^*}$ 
implies ${\bf G}(\partial {\bf e})={\bf G}(\partial {\bf e^*}),$
}
\mbox{thus ${\bf G}(\partial {\bf e})=S(\partial e_0)$ 
implies 
${\bf G}(\partial {\bf e^*})=S(\partial e_0)$})
\\
&{=}&
\Ex_{\varphi\in \Phi(mh)}
\Ex_{{\bf e^*}\in {\bf \Omega}_I
}
[
\brkt{
\Ex_{{\bf e}\in {\bf \Omega}_I
}[F_{e_0}\brkt{\bf e}
\llbracket
{\bf G}(\partial {\bf e})=S(\partial e_0)\rrbracket
|
{\bf e}\stackrel{
\partial {\bf G}/\varphi}{\approx}{\bf e^*}]
}^2
|{\bf G}(\partial {\bf e^*})=S(\partial e_0)
]
\\
&& 
\cdot
\brkt{
\prod_{e\in {\mathbb V}_{(k-1)}(S):e\subset e_0}
{\bf d}_{\bf G}^{(\delta)}(S\ang{e})
}
\brkt{
\prod_{e\in {\mathbb V}_{(k-1)}(S)}
{\bf d}_{\bf G}^{(\delta)}(S\ang{e})
}
\brkt{
\brkt{
\prod_{e\in {\mathbb V}_{(k-1)}(S),e\not\subset e_0}
{\bf d}_{\bf G}^{(\delta)}(S\ang{e})
}
+
{
1 \over m}
}
\\
&&
/
\brkt{
\prod_{e\in {\mathbb V}_{(k-1)}(S)}
{\bf d}_{\bf G}^{(-\delta)}(S\ang{e})
}^2
\\
&\le &
\Ex_{\varphi\in \Phi(mh)}
\Ex_{{\bf e^*}\in {\bf \Omega}_I
}
[
\brkt{
\Ex_{{\bf e}\in {\bf \Omega}_I
}[F_{e_0}\brkt{\bf e}
|
{\bf e}\stackrel{
\partial {\bf G}/\varphi}{\approx}{\bf e^*}]
}^2
|{\bf G}(\partial {\bf e^*})=S(\partial e_0)
]
\\
&& 
\cdot
\brkt{
1
+
{
1 \over m}
\brkt{
\prod_{e\in {\mathbb V}_{(k-1)}(S),e\not\subset e_0}
{\bf d}_{\bf G}^{(\delta)}(S\ang{e})
}^{-1}
}
\cdot
\brkt{
\prod_{e\in {\mathbb V}_{(k-1)}(S)}
{
{\bf d}_{\bf G}(S\ang{e})+
\delta(S\ang{e})
\over 
{\bf d}_{\bf G}(S\ang{e})
-
\delta(S\ang{e})
}
}^2.
\nonumber
\end{eqnarray*}
The assumption (\ref{070907b}) completes the proof of (\ref{070907c}).
\lemmaqed
\subsection{The body of our proof}
\begin{df}[Notation for this subsection]
Write $c_i({\bf G}):=\max_{I\in {{\mathfrak r}\choose i}}|{\rm C}_I({\bf G})|$ 
for $i\in [k]$.
For $\vec{b}=(b_i)_{i\in [k]}$ and an integer $m$,
we write 
$\vec{B}(\vec{b},m):=({B}_i(\vec{b},m))_{i\in [k]}$
where $
B_i(\vec{b},m):=
\prod_{j\in [0,k-i]}b_{i+j}^{{r-i\choose j}m^j}.
$ 
\end{df}
Recall (\ref{081230}). The $(k-1)$-regularization ${\bf G}/\varphi$ is
the $k$-bounded graph on ${\bf \Omega}
$ obtained from ${\bf G}$ by redefining the color of 
each edge ${\bf e}\in {\bf \Omega}_I$ with $I\in {\mathfrak{r}\choose i}$
by the 
$\brkt{\sum_{j=0}^{k-i}
{r-i\choose j}
m^{j}}$-dimensional 
vector
\begin{eqnarray*}
\brkt{{\bf G}/\varphi}({\bf e}):=({\bf G}({\bf e}\dot{\cup}{\bf f})
|
J\in {\mathfrak{r}\setminus I\choose [0,k-i]},
{\bf f}\in {\bf \Omega}_J \, \mbox{\rm with } \, 
{\bf f}\subset \varphi(\mathbb{D})
).
\end{eqnarray*}
Thus obviously if ${\bf G}$ is a $k$-bounded $\vec{b}$-colored graph then
\begin{eqnarray}
c_i({\bf G}/\varphi)
\le
B_i(\vec{b},m),\quad \forall i\in [k],\, \forall
\varphi\in\Phi(m).\label{lj24k}
\end{eqnarray}
(For example, $B_k(\vec{b},m)=b_k$ and $B_{k-1}(\vec{b},m)=b_{k-1}b_k^{(r-1)m}.$)
\par
Fix $0<\epsilon<1$ and $\vec{b}.$ We proceed by induction on $k.$
When $k=1,$ it is trivial as the remark after Theorem \ref{070908b}.
Let $k\ge 2.$
\par
$\bullet$
{\bf [Definition of the sample-size functions] }
Let
${m}^{(k-1)}_{k,h,\vec{b},\epsilon}(0)
:=0
$ and
${m}^{(i)}_{k,h,\vec{b},\epsilon}(n_i,\cdots,n_{k-2},0)
:={m}^{(i)}_{k-1,h,\vec{b},\epsilon}(n_i,\cdots,n_{k-2})
, \forall i\in [k-2]
$,
which is defined by the induction hypothesis
on $k-1$ of the theorem.
Define $\tilde{n}_{k,h,\vec{b},\epsilon}^{(k-1)}=\tilde{n}^{(k-1)}$
to be large enough so that
\begin{eqnarray}
Cb_k
\sqrt{
b_k\over \tilde{n}^{(k-1)}
}\le {\epsilon\over 4{r\choose k}}
\label{lj24g}
\end{eqnarray}
where
\begin{eqnarray}
C:=\sqrt{2}
{r\choose k}h^k
\brkt{b_k\over 2\sqrt{\epsilon_1}
}^{{r\choose k}h^k-1}
3^{
\sum_{j\in [k-1]}{r\choose j}h^j
}
\mbox{ and }
\epsilon_1:=\brkt{\epsilon\over 12\cdot 2^{k}b_k 
{r\choose k}
}^{2}.\label{lj24b}
\end{eqnarray}
(These expressions will appear in (\ref{lj25a}) and (\ref{lj24c}).)
Also let 
$\tilde{n}_{k,h,\vec{b},\epsilon}^{(j)}(n_{j+1},\cdots,n_{k-2},0)
:=\tilde{n}_{k-1,h,\vec{b},\epsilon}^{(j)}(n_{j+1},\cdots,n_{k-2})
$ for all $j\in [k-2]$.
\par
Given $n_{k-1}\ge 0$, we will inductively define functions
${m}^{(i)}_{k,h,\vec{b},\epsilon}(\bullet,\cdots,\bullet,n_{k-1}+1),
\forall i\in [k-1], $ and functions
$\tilde{n}^{(j)}_{k,h,\vec{b},\epsilon}(\bullet,\cdots,\bullet,n_{k-1}+1)
,\forall j\in [k-2], $ by using
${m}^{(\bullet)}_{k,h,\vec{b},\epsilon}({\bullet,\cdots,\bullet},n_{k-1}),$
 $m^{(\bullet)}_{k-1,\bullet,\bullet,\bullet}(\bullet,\cdots,\bullet),$
$\tilde{n}_{k,h,\vec{b},\epsilon}^{(\bullet)}(\bullet,\cdots,\bullet,n_{k-1}),$
and
$\tilde{n}_{k-1,\bullet,\bullet,\bullet}^{(\bullet)}(\bullet,\cdots,\bullet)$, 
as follows.
Let
\begin{eqnarray}
m&:=&
\prod_{i\in [k-1]}
\brkt{
B_i
\brkt{\vec{b},\sum_{j=1}^{k-1}
m_{k,h,\vec{b},\epsilon}^{(j)}
(\bar{n}^{(j)},\cdots,\bar{n}^{(k-2)},
n_{k-1}
)}
\over \sqrt{\epsilon_1}}^{{r\choose i}h^i}\label{lj24h}
\end{eqnarray}
where $\bar{n}^{(k-2)}:=\tilde{n}^{(k-2)}(n_{k-1})=
\tilde{n}^{(k-2)}_{k,h,\vec{b},\epsilon}(n_{k-1});$ $
\bar{n}^{(k-3)}:=\tilde{n}^{(k-3)}(\bar{n}^{(k-2)},n_{k-1});$ $
\cdots;$ $\bar{n}^{(j)}:=\tilde{n}^{(j)}( \bar{n}^{(j+1)}
,\cdots,\bar{n}^{(k-2)},n_{k-1}). $ (We will use the form
(\ref{lj24h}) only once in (\ref{lj22b}).) 
Define
$m_{k,h,\vec{b},\epsilon}^{(k-1)}(n_{k-1}+1)$ so that
\begin{eqnarray}
m_{k,h,\vec{b},\epsilon}^{(k-1)}(n_{k-1}+1)&:=&
\sum_{j=1}^{k-1}
m_{k,h,\vec{b},\epsilon}^{(j)}(\bar{n}^{(j)},\cdots,\bar{n}^{(k-2)},
n_{k-1}
)
+mh.
\label{lj24a}
\end{eqnarray}
Next, we define the remaining $k-2$ functions so that
\begin{eqnarray}
m_{k,h,\vec{b},\epsilon}^{(i)}(n_i,\cdots,n_{k-2},n_{k-1}+1)
&:=&m_{k-1,2h,\vec{b}^{*},\epsilon_1}^{(i)}(n_i,\cdots,n_{k-2}
), \quad \forall i\in [k-2] \label{lj24f}
\end{eqnarray}
where $\vec{b}^{*}=(b^{*}_i)_{i\in [k-1]}$ with
$
b^{*}_i
:=
B_i\brkt{\vec{b},m_{k,h,\vec{b},\epsilon}^{(k-1)}(n_{k-1}+1)
}.$
\  Finally we define
\begin{eqnarray}
\tilde{n}_{k,h,\vec{b},\epsilon}^{(j)}(n_{j+1},\cdots,n_{k-2},n_{k-1}+1)
:=
\tilde{n}_{k-1,2h,\vec{b}^{*},\epsilon_1}^{(j)}(n_{j+1},\cdots,n_{k-2}),
\quad \forall j\in [k-2].\label{lj24e}
\end{eqnarray}
(It will be easily seen that the three equalities 
$:=$ in (\ref{lj24a}),(\ref{lj24f}) and (\ref{lj24e}) 
can be replaced by $\ge$.)
\par
$\bullet$ {\bf [Definition of the error function]} For
$\vec{n}=(n^{(1)},\cdots,n^{(k-1)})$ and for
$\vec{\varphi}\in\Phi((m_{k,h,\vec{b},\epsilon}^{(j)}(n^{(j)}))_{j\in
[k-1]})$, we write ${\bf G^*}:={\bf G}/\vec{\varphi}$ and we
define a $(k,h)$-error function $
\delta=\delta_{k,h,\epsilon,{\bf G^*}}$ inductively as follows.
\par
Since (\ref{lj24k}) implies $c_i({\bf G}/\varphi_{k-1})
\le B_i(\vec{b}, m_{k,h,\vec{b},\epsilon}^{(k-1)}(n^{(k-1)}))$ and 
${\bf G^*}=({\bf G}/\varphi_{k-1})/(\varphi_i)_{i\in [k-2]}$, 
we apply 
the induction hypothesis on $k$ 
with (\ref{lj24f}) and (\ref{lj24e}) for ${\bf G}/\varphi_{k-1}$ 
and see 
that for the 
$\epsilon_1>0$ of (\ref{lj24b}), 
\begin{eqnarray}
\Ex_{\vec{n'}=(n^{(1)},\cdots,n^{(k-2)})}\Ex_{\vec{\varphi'}=(\varphi_i)_{i\in [k-2]}}
\left[
{\bf reg}_{k-1,2h}({\bf G^*})
\right]
\le 
\epsilon_1.
\nonumber
\end{eqnarray}
\\ 
Thus, there exists a function 
$\delta=\delta_{k-1,2h,\epsilon_1,{\bf G^*}}: {\rm
TC}({\bf G^*})=
\bigcup_{i<k}{\rm TC}_i({\bf G^*})\to [0,\infty)$ with the two property that (i) 
 for any $S\in {\mathcal S}_{k-1,2h,{\bf G^*}}$,
\begin{eqnarray}
\Prob_{\phi\in\Phi(h)}[{\bf G^*}(\phi(e))=S(e),\,\forall e\in
{\mathbb V}_{(k-1)}(S)] =\prod_{e\in {\mathbb V}_{(k-1)}(S)}{\bf
d}^{(\delta)}_{\bf G^*}(S\ang{e}) \label{lj25b}
\end{eqnarray}
and that (ii) for each fixed $\varphi_{k-1}\in
\Phi(m^{(k-1)}_{k,h,\vec{b},\epsilon}(n^{(k-1)}))$
\begin{eqnarray}
\Ex_{\vec{n'}=(n^{(1)},\cdots,n^{(k-2)})}\Ex_{\vec{\varphi'}=(\varphi_i)_{i\in [k-2]}}
\left[
\max_{I\in {\mathfrak{r}\choose [k-1]}}
|{\rm C}_I({\bf G^*})|
\Ex_{{\bf e}\in {\bf \Omega}_I}[
\delta({\bf G^*}\ang{\bf e})]
\right]
\le
\epsilon_1
\stackrel{(\ref{lj24b})}{<}
{\epsilon/2
}.
\label{lj24d}
\end{eqnarray}
(This $\delta_{k-1,2h,\epsilon_1,{\bf G^*}}$ depends (not only on $\varphi_{k-1}$ 
but also) 
on $\vec{n'}$ and $\vec{\varphi'}$.)
Define $\delta_{k,h,\epsilon,{\bf G^*}}(\vec{\mathfrak{c}}):=
\delta_{k-1,2h,\epsilon_1,{\bf G^*}}(\vec{\mathfrak{c}})$
for any $\vec{\mathfrak{c}}\in {\rm TC}_I({\bf G^*}), 
I\in {\mathfrak{r}\choose [k-1]}.$
\par
Before defining $\delta(\vec{\mathfrak{c}})$ for 
$\vec{\mathfrak{c}}\in {\rm TC}_k({\bf G^*})$,
we define \lq bad colors\rq\ ${\rm BAD}\subset {\rm TC}({\bf G^*}).$
For $I\in {\mathfrak{r}\choose [k]}$, we define ${\rm BAD}_I$ by the relation that
$\vec{\mathfrak{c}}=(\mathfrak{c}_J)_{J\subset I}\in {\rm BAD}_I$
if and only if
\begin{eqnarray}
\begin{array}{cccl}
 \delta( (\mathfrak{c}_J)_{J\subset I'}
)&\ge&
 \sqrt{\epsilon_1}/|{\rm C}_{I'}({\bf G^*})| & \mbox{ for some }
I'\subsetneq I,
 \mbox{ or }
\\
{\bf d}_{\bf G^*}((\mathfrak{c}_J)_{J\subset I^*}
)&\le &
2\sqrt{\epsilon_1}/|{\rm C}_{I^*}({\bf G^*})|& \mbox{ for some }
I^*\subset I.
\end{array}\label{lj24m}
\end{eqnarray}
Define ${\rm BAD}:=\bigcup_{I\in {\mathfrak{r}\choose [k]}}{\rm BAD}_I$.
\par
For $ \vec{\mathfrak{c}}=(\mathfrak{c}_J)_{J\subset I} \in {\rm
TC}_k({\bf G^*})$, we define, using $m$ and $C$ of
(\ref{lj24b}) and (\ref{lj24h}),
\begin{eqnarray}
\eta_{k,h}(\vec{\mathfrak{c}})&:=&
\Ex_{\varphi\in\Phi(mh)}\Ex_{{\bf e^*}\in {\bf \Omega}_I} [
\brkt{\Prob_{{\bf e}\in {\bf \Omega}_I}[ {\bf G^*}({\bf
e})=\mathfrak{c}_I | {\bf e}\stackrel{\partial {\bf
G^*}/\varphi}{\approx}{\bf e^*} ] - {\bf d}_{\bf
G^*}(\vec{\mathfrak c})}^2 |\, {\bf G^*}(\partial {\bf e^*})=
(\mathfrak{c}_J)_{J\subsetneq I}
],\label{061220}
\\
\delta_{k,h}(\vec{\mathfrak{c}})&:=&
\left\{
\begin{array}{cl}
1 & \mbox{ if } \vec{\mathfrak{c}}\in{\rm BAD}_I
,
\\
C\sqrt{\eta_{k,h}(\vec{\mathfrak{c}})},
& \mbox{otherwise.}\label{lj25}
\end{array}
\right.
\end{eqnarray}
\\
$\bullet$
{\bf [The qualification as an error function]}
Because of (\ref{lj25b}) and (\ref{lj25}),
it is enough for the first requirement (\ref{lj19}) 
to show that
\begin{eqnarray}
 \Prob_{\phi\in\Phi(h)}
[{\bf G^*}(\phi(e))=S(e)\,, \forall e\in {\mathbb V}(S)]=
\displaystyle\prod_{e\in {\mathbb V}(S)} \brkt{ 
{\bf d}_{\bf
G^*}(S\ang{e}) \dot{\pm} \delta(S\ang{e}) 
}
\label{090406}
\end{eqnarray}
or
\begin{eqnarray}
\Prob_{\phi\in\Phi(h)}[{\bf G^*}(\phi(e))=S(e),\,\forall e\in
{\mathbb V}_k(S)| {\bf G^*}(\phi(e))=S(e),\,\forall e\in {\mathbb
V}_{(k-1)}(S) ] =\prod_{e\in {\mathbb V}_k(S)}{\bf
d}^{(\delta)}_{\bf G^*}(S\ang{e})
\label{lj24}
\end{eqnarray}
for any $S\in {\mathcal S}_{k,h,{\bf G^*}}$.
Furthermore without loss of generality, we can assume the property that
\begin{eqnarray}
S\ang{e}\not\in{\rm BAD}
\mbox{ for any }
e\in {\mathbb V}(S). \label{090321}
\end{eqnarray}
(Indeed, we can show the case of (\ref{090321}) suffices by the induction 
on the number of bad edges in $S$.
Let a complex $S$ be given where $S$ contains a bad edge $e^*$.
Without loss of generality, assume that any visible edge $e\in \mathbb{V}(S)$ 
is not bad if $|e|<|e^*|.$
We construct a new complex $S^*$ from $S$ by recoloring 
all (bad) edges containing $e^*$ in the invisible color.
By the induction hypothesis, (\ref{090406}) holds for $S^*.$ 
Equality (\ref{090406}) means that 
the real number the left hand side suggests belongs to the 
interval which the right-hand side suggests. Denote by $[p^-, p^+]$ this interval.
Again we reconstruct $S$ from $S^*$ by recoloring some invisible edges in 
\lq original\rq\ 
bad colors.
By this process from $S^*$ to $S$, 
the left hand side of (\ref{090406}) will not increase 
(probably decrease because of added visible edges $e$)
and 
 the right-hand side will suggest interval $[0,p^+]$ because, for bad edges $e$, 
 ${\bf d}_{\bf
G^*}(S\ang{e}) \dot{\pm} \delta(S\ang{e}) =[0,1]$
by (\ref{lj25}). 
Then (\ref{090406}) holds not only for $S^*$ but also for $S$.)
\par
Fix such an $S\in {\mathcal S}_{k,h,{\bf G^*}}$. 
For any $e\in {\mathbb V}_J(S), J\subset \mathfrak{r}$, it follows from
(\ref{090321}) and (\ref{lj24m}) that
\begin{eqnarray}
{\bf d}^{(-\delta)}_{\bf G^*}(S\ang{e})
> {\sqrt{\epsilon_1}\over
|{\rm C}_J({\bf G^*})|}>0 
(\mbox{if }|J|< k)
\mbox{ and }
\delta(S\ang{e})
< {1\over 2}
{\bf d}_{\bf G^*}(S\ang{e})
 (\mbox{if }|J|< k).\label{lj22}
\end{eqnarray}
Clearly, 
$c_i({\bf G^*})=
c_i({\bf G}/\vec{\varphi})
\le c_i({\bf G}/^{k-1}(\varphi_{k-1}\dot{\cup}\cdots\dot{\cup}\varphi_1))$ and 
$|\mathbb{V}_i(S)|\le {r\choose i}h^i$.
Thus, it follows from (\ref{lj24h}) and (\ref{lj24k})  that
\begin{eqnarray}
{1\over m}
\le
\prod_{i\in [k-1]}
\brkt{\sqrt{\epsilon_1}\over c_i({\bf G^*})
}^{|{\mathbb V}_i(S)|}
\le^{(\ref{lj22})}
\prod_{i\in [k-1]}
\prod_{
e\in {\mathbb V}_{i}(S)
}
{\bf d}_{\bf G^*}^{(-\delta)}
(S\ang{e})
\le 
\prod_{e\in {\mathbb V}_{(k-1)}(S), e\not\subset e_0
}
{\bf d}_{\bf G^*}^{(-\delta)}
(S\ang{e})
\label{lj22b}
\end{eqnarray}
for any $e_0\in\mathbb{V}_k(S).$ 
Let $
F_e({\bf e}):=
\dbrkt{{\bf G^*}({\bf e})=S(e)
}-{\bf d}_{\bf G^*}(S\ang{e}).$ 
For any $\emptyset\not=D\subset {\mathbb V}_k(S)$,
we apply Lemma \ref{alj08} 
(where ${\bf G}:={\bf G^*}$) 
with any $S'\in {\mathcal S}_{k,h,{\bf G^*}}$ with
$\mathbb{V}_k(S')=D$ and $\mathbb{V}_{(k-1)}(S')=\mathbb{V}_{(k-1)}(S)$,
and see that
\begin{eqnarray}
&&
\brkt{
\Ex_{\phi\in\Phi(h)}\left[\left.
\prod_{e\in D}
\brkt{
\dbrkt{{\bf G^*}(\phi(e))=S(e)
}
-{\bf d}_{\bf G^*}(S\ang{e})
}
\right|
{\bf G^*}(\phi(e))=S(e)\forall e\in {\mathbb V}_{(k-1)}(S)
\right]}^2
\nonumber\\
&= &
\brkt{
\Ex_{\phi\in\Phi(h)}\left[\left.
\prod_{e\in D}
F_e(\phi(e))
\right|
{\bf G^*}(\phi(e))=S(e)\forall e\in {\mathbb V}_{(k-1)}(S)
\right]
}^2 
\nonumber
\\
&\le^{\rm Lem. \ref{alj08}}_{(\ref{lj22b}), (\ref{lj22})}
&
 \min_{e_0\in D}
2\cdot 3^{2|{\mathbb V}_{(k-1)}(S)|}
\cdot
\Ex_{\varphi\in \Phi(mh)} \Ex_{{\bf e^*}\in {\bf \Omega}_I } [
\brkt{ \Ex_{{\bf e}\in {\bf \Omega}_I }[F_{e_0}\brkt{\bf e}
|\,
 {\bf e}\stackrel{
\partial {\bf G^*}/\varphi}{\approx}{\bf e^*}]
}^2 |\,{\bf G^*}(\partial {\bf e^*})=S(\partial e_0) ] 
\nonumber\\
&\le^{(\ref{061220})} &
2\cdot 3^{2|{\mathbb V}_{(k-1)}(S)|}
\max_{e_0\in {\mathbb V}_k(S)}
\eta_{k,h}(S\ang{e_0})
.\label{lj23j}
\end{eqnarray}
Take an edge $e_0\in {\mathbb V}_k(S)$ which maximizes 
$\eta_{k,h}(S\ang{e_0}).$ Then it follows from 
Lemma \ref{lj23c} that 
\begin{eqnarray}
&& \Prob_{\phi\in\Phi(h)}[ {\bf G^*}(\phi(e))=S(e),\,\forall e\in
{\mathbb V}_k(S) | {\bf G^*}(\phi(e))=S(e)\forall e\in {\mathbb
V}_{(k-1)}(S) ]
\nonumber\\
&\stackrel{(\ref{lj23j})}{=}&
\prod_{e\in {\mathbb V}_k(S)}{\bf d}_{\bf G^*}(S\ang{e})
\dot{\pm} \sqrt{2}|{\mathbb V}_k(S)|
3^{|{\mathbb V}_{(k-1)}(S)|}
\sqrt{\eta_{k,h}(S\ang{e_0})}
\quad (\mbox{by taking }D:=\mathbb{V}_k(S))
\nonumber\\
&=^{(\ref{090321})}_{(\ref{lj24m})}&
\brkt{
{\bf d}_{\bf G^*}(S\ang{e_0})
\dot{\pm} {\sqrt{2}
{r\choose k}h^k
3^{\sum_{j\in [k-1]}{r\choose j}h^j
}
\sqrt{\eta_{k,h}(S\ang{e_0})}
\over
\brkt{2\sqrt{\epsilon_1}/c_k({\bf G^*})}^{|{\mathbb V}_k(S)|-1}
}
}
\prod_{e\in {\mathbb V}_k(S),e\not=e_0}
{\bf d}_{\bf G^*}(S\ang{e})
\nonumber\\
&=^{(\ref{lj24b})}_{(\ref{lj25})}
&\prod_{e\in {\mathbb V}_k(S)}
({\bf d}_{\bf G^*}(S\ang{e})
\dot{\pm}\delta(S\ang{e}))
\label{lj25a}
\end{eqnarray}
where for the last equality we use the fact that $c_k({\bf G^*})=c_k({\bf G})\le b_k$
(cf. (\ref{lj24k})).
\bigskip
\\
$\bullet$
{\bf [Bounding the average error size] }
With the abbreviation $a_n:=m^{(k-1)}_{k,h,\vec{b},\epsilon}(n)$,
for any $I\in {\mathfrak{r}\choose k}$, the linearity of expectation gives us that 
\begin{eqnarray}
&&
\brkt{
\Ex_{\vec{n},\vec{\varphi}}
\Ex_{{\bf \tilde{e}}\in {\bf \Omega}_I}[
\sqrt{
\eta_{k,h}({\bf G^*}\ang{\bf \tilde{e}})}
]
}^2
\nonumber
\\
&\le &
\Ex_{\vec{n},\vec{\varphi}}
\Ex_{{\bf \tilde{e}}\in {\bf \Omega}_I}[
\eta_{k,h}({\bf G^*}\ang{\bf \tilde{e}})
]\quad\quad (\mbox{by Cauchy-Schwarz or }  \Ex[X]^2\le \Ex[X^2])
\nonumber
\\
&\stackrel{(\ref{061220})}{=}
& \Ex_{\vec{n},\vec{\varphi}, {\bf \tilde{e}}}
\Ex_{\varphi\in\Phi(mh)}\Ex_{{\bf e^*}\in {\bf \Omega}_I} [
\brkt{\Prob_{{\bf e}\in {\bf \Omega}_I}[ {\bf G^*}({\bf e})={\bf
G^*}({\bf \tilde{e}}) |\, {\bf e}\stackrel{\partial {\bf
G^*}/\varphi}{\approx}{\bf e^*} ] - {\bf d}_{\bf G^*}({\bf
G^*}\ang{\bf \tilde{e}} )}^2 |\,
 {\bf e^*}
\stackrel{\partial {\bf G^*}}{\approx}{\bf \tilde{e}}
]
\nonumber
\\
&\le &\Ex_{\vec{n},\vec{\varphi},
{\bf \tilde{e}}}
\sum_{\mathfrak{c}_I\in {\rm C}_I({\bf G^*})}
\Ex_{\varphi,
{\bf e^*}}
[
\brkt{\Prob_{{\bf e}\in {\bf \Omega}_I}[
{\bf G^*}({\bf e})=\mathfrak{c}_I
|\, {\bf e}\stackrel{\partial {\bf G^*}/\varphi}{\approx}{\bf e^*}
] - \Prob_{{\bf e}\in {\bf \Omega}_I}[ {\bf G^*}({\bf
e})=\mathfrak{c}_I
|\, {\bf e}\stackrel{\partial {\bf G^*}}{\approx}{\bf \tilde{e}} ]
}^2 |\,
 {\bf e^*}
\stackrel{\partial {\bf G^*}}{\approx}{\bf \tilde{e}}
]
\nonumber
\\
&= & \sum_{\mathfrak{c}_I\in {\rm C}_I({\bf G})}
\Ex_{\vec{n},\vec{\varphi}, {\bf \tilde{e}}} \left[ 
\Ex_{\varphi,
{\bf e^*}} [ \brkt{\Prob_{{\bf e}\in {\bf \Omega}_I}[ {\bf G}({\bf
e})=\mathfrak{c}_I |\, {\bf e}\stackrel{\partial {\bf
G^*}/\varphi}{\approx}{\bf e^*} ]}^2 |
 {\bf e^*}
\stackrel{\partial {\bf G^*}}{\approx}{\bf \tilde{e}} ] 
+\brkt{
\Prob_{{\bf e}\in {\bf \Omega}_I}[ {\bf G}({\bf e})=\mathfrak{c}_I
|\, {\bf e}\stackrel{\partial {\bf G^*}}{\approx}{\bf \tilde{e}} ]}^2
\right. 
\nonumber\\
&&\left.
-2\Ex_{\varphi,
{\bf e^*}} [ 
\Prob_{{\bf e}\in {\bf \Omega}_I}[ {\bf G}({\bf
e})=\mathfrak{c}_I |\, {\bf e}\stackrel{\partial {\bf
G^*}/\varphi}{\approx}{\bf e^*} ]
|
 {\bf e^*}
\stackrel{\partial {\bf G^*}}{\approx}{\bf \tilde{e}}
]\cdot
\Prob_{{\bf e}\in {\bf \Omega}_I}[ {\bf G}({\bf e})=\mathfrak{c}_I
|\, {\bf e}\stackrel{\partial {\bf G^*}}{\approx}{\bf \tilde{e}} ]
\right] \nonumber
\\
&= & \sum_{\mathfrak{c}_I\in {\rm C}_I({\bf G})}
\Ex_{\vec{n},\vec{\varphi}, {\bf \tilde{e}}} \left[ \Ex_{\varphi,
{\bf e^*}} [ \brkt{\Prob_{{\bf e}\in {\bf \Omega}_I}[ {\bf G}({\bf
e})=\mathfrak{c}_I |\, {\bf e}\stackrel{\partial {\bf
G^*}/\varphi}{\approx}{\bf e^*} ]}^2 |
 {\bf e^*}
\stackrel{\partial {\bf G^*}}{\approx}{\bf \tilde{e}} ] -\brkt{
\Prob_{{\bf e}\in {\bf \Omega}_I}[ {\bf G}({\bf e})=\mathfrak{c}_I
|\, {\bf e}\stackrel{\partial {\bf G^*}}{\approx}{\bf \tilde{e}} ]
)}^2 \right]
 \nonumber
\\
&= &|{\rm C}_I({\bf G})|
 \Ex_{\mathfrak{c}_I\in {\rm C}_I({\bf G})} \Ex_{\vec{n},\vec{\varphi}, {\bf
\tilde{e}}} \left[ 
\Ex_{\varphi\in\Phi(mh)} [ \brkt{\Prob_{{\bf e}}[ {\bf G}({\bf e})=\mathfrak{c}_I |\, {\bf e}\stackrel{\partial
({\bf G}/\vec{\varphi}) /\varphi}{\approx}{\bf \tilde{e}} ]}^2 ]
-\brkt{ \Prob_{{\bf e}}[ {\bf G}({\bf e})=\mathfrak{c}_I | {\bf
e}\stackrel{\partial ({\bf G}/\vec{\varphi}) }{\approx}{\bf
\tilde{e}} ] )}^2
 \right] 
\nonumber
\\
&\stackrel{(*)}{\le} &b_k \Ex_{0\le n<\tilde{n}^{(k-1)}} \Ex_{{\bf
\tilde{e}},\mathfrak{c}_I}
\left[ 
\Ex_{\phi'\in \Phi(a_{n+1})}[ \brkt{
\Prob_{{\bf e}}[ {\bf G}({\bf e})=\mathfrak{c}_I |\, {\bf
e}\stackrel{\partial {\bf G}/\phi'}{\approx}{\bf \tilde{e}} ] }^2
] - \Ex_{{\phi}\in\Phi(a_n)}[ \brkt{ \Prob_{{\bf e}}[ {\bf G}({\bf
e})=\mathfrak{c}_I |\, {\bf e}\stackrel{\partial {\bf
G}/{\phi}}{\approx}{\bf \tilde{e}} ] }^2 ]
\right]
\nonumber\\
&= & {b_k\over \tilde{n}^{(k-1)}} \sum_{n=0}^{\tilde{n}^{(k-1)}-1}
\Ex_{{\bf \tilde{e}},\mathfrak{c}_I}
\left[ \Ex_{\phi\in \Phi(a_{n+1})}[
\brkt{ \Prob_{{\bf e}}[ {\bf G}({\bf e})=\mathfrak{c}_I |\,
 {\bf
e}\stackrel{\partial {\bf G}/\phi}{\approx}{\bf \tilde{e}} ] }^2 ]
- \Ex_{{\phi}\in\Phi(a_n)}[ \brkt{ \Prob_{{\bf e}}[ {\bf G}({\bf
e})=\mathfrak{c}_I | {\bf e}\stackrel{\partial {\bf
G}/{\phi}}{\approx}{\bf \tilde{e}} ] }^2 ]\right]
\nonumber\\
&= & {b_k\over \tilde{n}^{(k-1)}} \Ex_{{\bf
\tilde{e}},\mathfrak{c}_I}
\left[ \Ex_{\phi\in
\Phi(a_{\tilde{n}^{(k-1)}})}[ \brkt{ \Prob_{{\bf e}}[ {\bf G}({\bf
e})=\mathfrak{c}_I |\,
 {\bf e}\stackrel{\partial {\bf
G}/\phi}{\approx}{\bf \tilde{e}} ] }^2 ] -
\Ex_{{\phi}\in\Phi(a_0)}[ \brkt{ \Prob_{{\bf e}}[ {\bf G}({\bf
e})=\mathfrak{c}_I |\, {\bf e}\stackrel{\partial {\bf
G}/{\phi}}{\approx}{\bf \tilde{e}} ] }^2 ]\right]
\nonumber\\
&\le&{b_k\over \tilde{n}^{(k-1)}},\label{lj20}
\end{eqnarray}
where in the above (*) 
we use the property that, after $n=n^{(k-1)}$ is chosen, it
follows from 
(\ref{lj24a}) that $a_{n+1}\ge
mh+\sum_{j=1}^{k-1}m^{(j)}(n^{(j)},\cdots,n^{(k-2)},
n^{(k-1)}) \ge a_n$ (for all possible 
$n^{(1)},\cdots,n^{(k-2)}$)
(cf. definition of $\bar{n}^{(k-2)},\cdots,\bar{n}^{(1)}$ just after 
(\ref{lj24h}) )
and
that if $\phi'({\mathbb D})\supset (\bigcup_{i\in
[k-1]}\varphi_i({\mathbb D})) \cup \varphi({\mathbb D})$ then
${\bf e}\stackrel{\partial {\bf G}/\phi'}{\approx}{\bf \tilde{e}}
$ implies ${\bf e}\stackrel{\partial ({\bf
G}/(\varphi_i)_i)/\varphi}{\approx}{\bf \tilde{e}} $ 
(thus $
\Ex_{\vec{\varphi}}
\Ex_{\varphi\in\Phi(mh)} 
[ \brkt{\Prob_{{\bf e}}[ {\bf G}({\bf e})=\mathfrak{c}_I |\, {\bf e}\stackrel{\partial
({\bf G}/\vec{\varphi}) /\varphi}{\approx}{\bf \tilde{e}} ]}^2 ]
\stackrel{(\ref{080827})}{\le}
\Ex_{\phi'\in \Phi(a_{n+1})}[ \brkt{
\Prob_{{\bf e}}[ {\bf G}({\bf e})=\mathfrak{c}_I |\, {\bf
e}\stackrel{\partial {\bf G}/\phi'}{\approx}{\bf \tilde{e}} ] }^2
] 
$)
and further,
that 
${\bf e}\stackrel{\partial ({\bf
G}/(\varphi_i)_i)}{\approx}{\bf \tilde{e}} $
 implies
${\bf e}\stackrel{\partial {\bf G}/\phi}{\approx}{\bf \tilde{e}}
$ where $\phi=\varphi_{k-1}$
(thus 
$
\Ex_{\vec{\varphi}}
[
\brkt{ \Prob_{{\bf e}}[ {\bf G}({\bf e})=\mathfrak{c}_I | {\bf
e}\stackrel{\partial ({\bf G}/\vec{\varphi}) }{\approx}{\bf
\tilde{e}} ] )}^2
]
\stackrel{(\ref{080827})}{\ge}
 \Ex_{{\phi}\in\Phi(a_n)}[ \brkt{ \Prob_{{\bf e}}[ {\bf G}({\bf
e})=\mathfrak{c}_I |\, {\bf e}\stackrel{\partial {\bf
G}/{\phi}}{\approx}{\bf \tilde{e}} ] }^2 ]
$)
.
\par
Thus, for any $I\in {\mathfrak{r}\choose k},$ we see that
\begin{eqnarray}
&& 
\Ex_{\vec{n},\vec{\varphi}}[ |{\rm C}_I({\bf
G}/\vec{\varphi})| \,\Ex_{{\bf e}\in {\bf
\Omega}_I}[\delta_{k,h}({\bf G^*}\ang{\bf e})]]
\nonumber
\\
&\stackrel{(\ref{lj25})}{\le} &
b_k\Ex_{\vec{n},\vec{\varphi}}\left[
{
C
\Ex_{
{\bf e}\in {\bf \Omega}_I
}[
\sqrt{\eta_{k,h}({\bf G^*}\ang{\bf e}
)}]
}
+1\cdot \Prob_{{\bf e}\in {\bf \Omega}_I}[{\bf G^*}\ang{\bf e}\in {\rm BAD}_I]
\right]
\nonumber
\\
&\le^{(\ref{lj20})}_{(\ref{lj24m})} &
b_k
\brkt{
{
C
\sqrt{
b_k\over \tilde{n}^{(k-1)}
}
}
+
\Ex_{\vec{n},\vec{\varphi}}
\left[
\sum_{J\subsetneq I}
\Prob_{{\bf e}\in {\bf \Omega}_J}[\delta({\bf G^*}\ang{\bf e})\ge
{\sqrt{\epsilon_1}\over |{\rm C}_J({\bf G^*})|}
]
+\sum_{J\subset I}
\Prob_{{\bf e}\in {\bf \Omega}_J}
[{\bf d}_{\bf G^*}({\bf G^*}\ang{\bf e})\le
{2\sqrt{\epsilon_1}\over |{\rm C}_J({\bf G^*})|}
]
\right]
}
\nonumber
\\
&\le^{(\ref{lj24d})}_{(**)} &
b_k\brkt{
C
\sqrt{
b_k\over \tilde{n}^{(k-1)}
}
+2^k
\sqrt{\epsilon_1}
+2^k\cdot 2
\sqrt{\epsilon_1}
}
\le_{(\ref{lj24b})}^{(\ref{lj24g})}
{\epsilon
\over 2{r\choose k}}.\label{lj24c}
\end{eqnarray}
where in the above (**) we use (\ref{lj26}) and the fact that 
\begin{eqnarray}
&&
\Prob_{{\bf e}\in {\bf \Omega}_J}\left[
\Prob_{{\bf e'}\in {\bf \Omega}_J}
[{\bf G^*}({\bf e'})
={\bf G^*}({\bf e})
|{\bf e'}\stackrel{\partial {\bf G^*}}{\approx}{\bf e}
]
\le
{2\sqrt{\epsilon_1}\over |{\rm C}_J({\bf G^*})|}
\right]
\nonumber
\\
&=&\sum_{\mathfrak{c}_J\in {\rm C}_J({\bf G^*})}
\Prob_{{\bf e}\in {\bf \Omega}_J}\left[
{\bf G^*}({\bf e})=\mathfrak{c}_J \mbox{ and }
\Prob_{{\bf e'}\in {\bf \Omega}_J}
[{\bf G^*}({\bf e'})
=\mathfrak{c}_J
|{\bf e'}\stackrel{\partial {\bf G^*}}{\approx}{\bf e}
]
\le
{2\sqrt{\epsilon_1}\over |{\rm C}_J({\bf G^*})|}
\right]\nonumber
\\
&\le &\sum_{\mathfrak{c}_J\in {\rm C}_J({\bf G^*})}
1\cdot
\Prob_{{\bf e}\in {\bf \Omega}_J}\left[
{\bf G^*}({\bf e})=\mathfrak{c}_J \left|
\Prob_{{\bf e'}\in {\bf \Omega}_J}
[{\bf G^*}({\bf e'})
=\mathfrak{c}_J
|{\bf e'}\stackrel{\partial {\bf G^*}}{\approx}{\bf e}
]
\le
{2\sqrt{\epsilon_1}\over |{\rm C}_J({\bf G^*})|}
\right.
\right]
\nonumber
\\
&= &\sum_{\mathfrak{c}_J\in {\rm C}_J({\bf G^*})}
\Ex_{{\bf e}\in {\bf \Omega}_J}\left[
\Prob_{{\bf \tilde{e}}\in {\bf \Omega}_J}
[
{\bf G^*}({\bf \tilde{e}})=\mathfrak{c}_J 
|
{\bf \tilde{e}}\stackrel{\partial {\bf G^*}}{\approx}
{\bf e}
]
\left|
\Prob_{{\bf e'}\in {\bf \Omega}_J}
[{\bf G^*}({\bf e'})
=\mathfrak{c}_J
|{\bf e'}\stackrel{\partial {\bf G^*}}{\approx}{\bf e}
]
\le
{2\sqrt{\epsilon_1}\over |{\rm C}_J({\bf G^*})|}
\right.
\right]\nonumber
\\
&& \hspace{1.5in}\mbox{($\because$ the conditional part depends only on ${\bf G^*}(\partial {\bf e})$)}
\nonumber
\\
&\le  &\sum_{\mathfrak{c}_J\in {\rm C}_J({\bf G^*})}
\Ex_{{\bf e}\in {\bf \Omega}_J}\left[
{2\sqrt{\epsilon_1}\over |{\rm C}_J({\bf G^*})|}
\left|
\Prob_{{\bf e'}\in {\bf \Omega}_J}
[{\bf G^*}({\bf e'})
=\mathfrak{c}_J
|{\bf e'}\stackrel{\partial {\bf G^*}}{\approx}{\bf e}
]
\le
{2\sqrt{\epsilon_1}\over |{\rm C}_J({\bf G^*})|}
\right.
\right]
=2\sqrt{\epsilon_1}.\label{070908}
\end{eqnarray}
Thus we obtain that 
\begin{eqnarray*}
&&
\Ex_{\vec{n},\vec{\varphi}}[ 
{\bf reg}_{k,h}({\bf G}/\vec{\varphi})
]\nonumber
\\
&\le &
\Ex_{\vec{n},\vec{\varphi}}[ \max_{I\in {\mathfrak{r}\choose [k]}}
|{\rm C}_{I}({\bf
G}/\vec{\varphi})| \,\Ex_{{\bf e}\in {\bf
\Omega}_I}[\delta({\bf G^*}\ang{\bf e})]]
\nonumber
\\
&\le &
\Ex_{\vec{n},\vec{\varphi}}[ 
\max_{I\in {\mathfrak{r}\choose [k-1]}}
|{\rm C}_{I}({\bf
G}/\vec{\varphi})| \,\Ex_{{\bf e}\in {\bf
\Omega}_I}[\delta({\bf G^*}\ang{\bf e})]]
+
\Ex_{\vec{n},\vec{\varphi}}[ \sum_{I\in {\mathfrak{r}\choose k}}
|{\rm C}_I({\bf
G}/\vec{\varphi})| \,\Ex_{{\bf e}\in {\bf
\Omega}_I}[\delta({\bf G^*}\ang{\bf e})]]
\nonumber
\\
&{\le}^{(\ref{lj24d})}_{(\ref{lj24c})} &
\epsilon/2
+\sum_{I\in {\mathfrak{r}\choose k}}{\epsilon\over 2
{r\choose k}}=\epsilon.
\end{eqnarray*}
It shows the second requirement (\ref{lj26}) for function $\delta$,
completing the proof of the main theorem.
\qed
\section{The Removal Lemma and Proof of Theorem \ref{24}}
While there had been known that
some strong versions of hypergraph regularity lemmas imply
Szemer\'edi's theorem (\cite{FR02}) before they were proven,
Solymosi \cite{So03,So04} inspired by Erd\H{o}s and Graham
showed that they also yield a combinatorial proof of
Theorem \ref{24}.
We will describe his argument for completeness and for seeing the
length of the entire proof of Theorem \ref{24}.
\begin{df}[$k$-uniform graphs]
A {\bf $k$-uniform $b_k$-colored ($\mathfrak{r}$-partite hyper)graph} is
a $k$-bounded $(b_i)_{i\in [k]}$-colored graph such that (1)
 if $i<k$ then $b_i=1$ and the unique color is called invisible
and (2) for each $I$ with $|I|=k$, there is at most one
 index-$I$ color which is called invisible. Denote by $\mathbb{V}(F)$ the set of
visible edges of a $k$-uniform graph $F$, where a visible edge
means an edge whose color is not invisible. Such a graph is called {\bf
$h$-vertex} if each partite set contains exactly $h$ vertices.
\end{df}
\begin{thm}[The Removal Lemma]\label{a01a}
For any $r\ge k,\, h,\,\vec{b}=(b_i)_{i\in [k]},$ and for any
$\epsilon>0$, there exists a constant
$c=c_{\ref{a01a}}(r,k,h,\vec{b},\epsilon)>0$ with the
following property.\par
Let ${\bf G}$ be a $k$-bounded $\vec{b}$-colored
($\mathfrak{r}$-partite hyper)graph on ${\bf \Omega}=({\bf
\Omega}_i)_{i\in \mathfrak{r}}.$ Let $F$ be an $h$-vertex
$k$-uniform $(b_k-1)$-colored ($\mathfrak{r}$-partite hyper)graph.
Then at least one of the following two holds.
\\
{\rm (i)}
 There exists a $k$-bounded $\vec{b}$-colored
($\mathfrak{r}$-partite hyper)graph ${\bf G'}$ on
${\bf \Omega}$ such that
\begin{eqnarray*}
\Prob_{{\bf e}\in {\bf \Omega}_I }[{\bf G'}({\bf e})\not={\bf
G}({\bf e})]\le \epsilon, \quad \forall I\in {\mathfrak{r}\choose
k}
\quad\mbox{
and }\quad
\Prob_{\phi\in\Phi(h)}[{\bf G'}(\phi(e))=F(e),\,\forall e\in
{\mathbb V}(F)]=0.
\end{eqnarray*}
{\rm (ii) }
$ \Prob_{\phi\in\Phi(h)}[{\bf G}(\phi(e))=F(e),\, \forall e\in
{\mathbb V}(F)]\ge c.
$
\end{thm}
\begin{proof}[Proof][Tool: Corollary \ref{lj21}]
 Let $\varepsilon\le \brkt{\epsilon \over 3\cdot 2^{k}}^{2}$, which is different 
 from $\epsilon.$
 Corollary \ref{lj21} gives
constants $\widetilde{m}_1,\cdots,\widetilde{m}_{k-1}$ such that,
given ${\bf G}$, there exist constants $m_1\le
\widetilde{m}_1,\cdots,m_{k-1}\le \widetilde{m}_{k-1}$ together
with $\vec{\varphi}\in\Phi(m_1,\cdots,m_{k-1})$ and with a
$(k,h)$-error function $\delta=\delta_{\vec{\varphi}}$ of ${\bf
G^*}:={\bf G}/\vec{\varphi}$ for which
\begin{eqnarray}
\Ex_{{\bf e}\in {\bf \Omega}_I}[\delta({\bf G^*}\ang{\bf e})]\le
\varepsilon/|{\rm C}_I({\bf G^*})|, \quad \forall I\in
{\mathfrak{r}\choose [k]}. \label{a01}
\end{eqnarray}
For $I\in {\mathfrak{r}\choose k},$ define
${\rm BAD}_I\subset {\rm TC}_I({\bf G^*})$
by the relation that $\vec{\mathfrak{c}}=(\mathfrak{c}_J)_{J\subset I}\in
{\rm BAD}_I$ if and only if
there exists an $I'\subset I$ such that
$
{\bf d}_{{\bf G^*}}(({\mathfrak c}_J)_{J\subset I'}
)\le 2\sqrt{\varepsilon}/|{\rm C}_{I'}({\bf G^*})| \mbox{ or that }
\delta(({\mathfrak c}_J)_{J\subset I'})\ge \sqrt{\varepsilon}/|{\rm C}_{I'}({\bf G^*})|.
$
For each $I\in {\mathfrak{r}\choose k},$ there exists a color
${\mathfrak c}^*_I\in {\rm C}_I({\bf G})\setminus {\rm C}_I(F) $
since  $F$ is $(b_k-1)$-colored. We replace each
$\vec{\mathfrak{c}}=(\mathfrak{c}_J)_{J\subset I}\in {\rm BAD}_I$
by $\vec{\mathfrak{c}^*}=(\mathfrak{c}^*_J)_{J\subset I}$ where
$\mathfrak{c}^*_J:=\mathfrak{c}_J$ for any $J\subsetneq I.$ Denote
the resulting graph by ${\bf G'}$. Then for each $I\in
{\mathfrak{r}\choose k}$ the same argument as  in 
(\ref{lj24c}) and (\ref{070908}) gives that 
\begin{eqnarray}
&&
\Prob_{{\bf e}\in {\bf \Omega}_I}[{\bf G'}({\bf e})\not={\bf G}({\bf e})]
\nonumber\\
&=&\Prob_{{\bf e}\in {\bf \Omega}_I}[
{\bf G}({\bf e})\in {\rm BAD}_I]\nonumber
\\
&\le& \sum_{I'\subset I}\brkt{
\Prob_{{\bf e}\in {\bf \Omega}_{I'}}[{\bf d}_{\bf G^*}({\bf G^*}\ang{\bf e})\le
{2\sqrt{\varepsilon}\over
|{\rm C}_{I'}({\bf G^*})|}]
+
\Prob_{{\bf e}\in {\bf \Omega}_{I'}}[\delta_{\bf G^*}({\bf G^*}\ang{\bf e})\ge
{\sqrt{\varepsilon}\over
|{\rm C}_{I'}({\bf G^*})|}]}\nonumber
\\
&\le^{(\ref{a01})}_{(\ref{lj24c},\ref{070908})} & 2^k\brkt{2\sqrt{\varepsilon}+
\sqrt{\varepsilon}}=3\cdot 2^{k}\sqrt{\varepsilon}
\le \epsilon.\label{a02}
\end{eqnarray}
Consider an $S\in {\mathcal S}_{k,h,{\bf G^*}}$ such that
$\mathbb{V}_k(S)=\mathbb{V}(F)$ and such that $S(e)=F(e)$ for all
$e\in \mathbb{V}_k(S).$ Denote by ${\mathcal S}^*$ the set of such
$S$ with the additional property that $S\ang{e}\not\in \bigcup_I
{\rm BAD}_I$ for any $e\in {\mathbb V}_k(S).$ Then our way of
recoloring gives that
\begin{eqnarray*}
\Prob_{\phi\in\Phi(h)}[{\bf G}(\phi(e))=F(e),\,\forall e\in
\mathbb{V}(F)] &\ge & \Prob_{\phi\in\Phi(h)}[{\bf
G'}(\phi(e))=F(e)\forall e\in \mathbb{V}(F)]
\\
&= & \sum_{S\in {\mathcal S}^*} \Prob_{\phi\in\Phi(h)}[{\bf
G^*}(\phi(e))=S(e),\,\forall e\in \mathbb{V}(S)]
\\
&=& \sum_{S\in {\mathcal S}^*} \prod_{e\in {\mathbb V}(S)}({\bf
d}_{\bf G^*}(S\ang{e})\dot{\pm}\,\delta(S\ang{e}))
\\
&\ge &
\sum_{S\in {\mathcal S}^*}\prod_{I\in {\mathfrak{r}\choose [k]}}
\prod_{e\in {\mathbb V}_I(S)}
{2\sqrt{\varepsilon}-\sqrt{\varepsilon}\over
|{\rm C}_I({\bf G^*})|}
\\
&\ge &
|{\mathcal S}^*|\prod_{I\in {\mathfrak{r}\choose [k]}}
\brkt{\sqrt{\varepsilon}\over
(b_{|I|}+\cdots+b_k)^{r(\widetilde{m}_{|I|}+\cdots+\widetilde{m}_{k-1})\choose k}
}^{h^{|I|}}
.
\end{eqnarray*}
Therefore if ${\mathcal S}^*=\emptyset$ then the first equality in the above 
with (\ref{a02}) gives the first condition.
Otherwise the second condition holds.
\end{proof}
For an integer $m$, we write $[m]_0:=[0,m-1]=\{0,1,\cdots,m-1\}.$
Write $E_r:=\{ (\underbrace{0,\cdots,0}_{i-1},1, 0,\cdots,0)\in
\integerset^{r}\,|\, i\in [r] \}.$
\begin{lemma}
\label{bnj08} For any $\delta>0$ and $k\ge 1,$ there exists an
$\epsilon>0$ satisfying the following. If an integer $N$ is
sufficiently large then for any subset $S\subset T(N,k):=
\{x=(x_0,\cdots,x_k)\in [N]_0^{k+1} | x_0+\cdots+x_k=N-1\} $ with
$|S|\ge \delta N^k$, there exists $a=(a_0,\cdots,a_k)\in
\integerset^{k+1}\setminus T(N,k)$
 with
$a+cE_{k+1}\subset S$ where
$c:
=N-1-\sum_{i=0}^k a_i\not=0$.
 Furthermore, there are at least $\epsilon N^{k+1}$ of such vectors $a.$
\end{lemma}
\begin{proof}[Proof][Tool: Theorem \ref{a01a}]
Let $S\subset T(N,k).$
Let $\mathfrak{r}:=\{0,\cdots,k\}$
and $({\bf \Omega}_i,{\mathcal B}_i,\mu_i):=([N]_0,2^{[N]_0},\mu_i(\bullet)=1/|\bullet|)$
for $i\in \mathfrak{r}.$
Define a $(b_1=1,\cdots,b_{k-1}=1,b_k=2)$-colored
$k$-bounded $\mathfrak{r}$-partite hypergraph ${\bf G}$ with vertex sets
${\bf \Omega}=({\bf \Omega}_i)_{i\in\mathfrak{r}}
$ so that
for each $I\in {\mathfrak{r}\choose k}$ and for each $k$-tuple
${\bf e}=(x_i\in [N]_0)_{i\in I}\in {\bf \Omega}_I$
, ${\bf e}$ is {\em red} if and only if
there exists
$v=(v_i)_{i\in [0,k]}\in S$ such that $v_i=x_i$ for any $i\in I.$
\par
Let $F$ be a $1$-vertex $k$-uniform $1$-colored graph on vertices
$V(F)=(V_i(F)=\{i\})_{i\in \mathfrak{r}}$ such that all the $k+1$
visible edges of $F$ are red. We say that
$\phi\in\Phi(1)=\Phi([0,k],{\bf \Omega})$ is {\em red} (in ${\bf
G}$) if and only if $ {\bf G}(\phi(e))={\rm red}$ for any $e\in
\mathbb{V}(F).$ We also say that a red $\phi\in\Phi(1)$ is {\em
degenerate} if and only if $(\phi(i))_{i\in [0,k]}\in S.$ 
Suppose that there exists a graph ${\bf G'}$ such that $\Prob_{{\bf e}\in {\bf
\Omega}_I}[{\bf G'}({\bf e})\not={\bf G}({\bf e})] \le 0.99
\delta/(k+1)$ for any $I\in {\mathfrak{r}\choose k}$ and
$\Prob_{\phi\in\Phi(1)}[\phi \mbox{ is red in } {\bf G'}]=0$.
Then  $|S|= |\{ \phi\in\Phi(1): \phi \mbox{
is degenerate}\}| \le \sum_{I\in {\mathfrak{r}\choose k}} |\{ {\bf
e}\in {\bf \Omega}_I : {\bf G'}({\bf e})\not={\bf G}({\bf e}) \}|
\le 0.99\delta N^k<|S|,$ where (in the first inequality) 
we use the fact that one cannot
delete two distinct degenerate $\phi$'s by recoloring one red edge in ${\bf G}$. 
Therefore, such a graph ${\bf G'}$ does not exist and 
Theorem \ref{a01a} gives a constant
$c^*=c_{\ref{a01a}}(r=k+1,k=k,h=1,\vec{b}=
(1,\cdots,1,2),\epsilon:=0.99\delta/(k+1))>0$ such that
\begin{eqnarray*}
\Prob_{\phi\in\Phi(1)}[
\phi \mbox{ is non-degenerate red}]
&\ge&
\Prob_{\phi\in\Phi(1)}[{\bf G}(\phi(e))=F(e)\forall e\in \mathbb{V}(F)]
-\Prob_{\phi}[
 \phi \mbox{ is degenerate}
]
\\
&\ge &
c^*-|S|/N^{k+1}\ge c^*-1/N.
\end{eqnarray*}
Thus, if $N\ge 1/0.9c^*$ then there exist $0.1c^* N^{k+1}$
non-degenerate red $\phi\in\Phi(1)$. 
Observe that a non-degenerate
red $\phi$ yields the desired $a+cE_{k+1}\subset S$ with
$a:=(\phi(i))_{i\in [0,k]}$ and $ c:=N-1-\sum_{i=0}^k
\phi(i)\not=0$ since if $c=0$ then it is degenerate.
\end{proof}
{\bf Proof of Theorem \ref{24}:} [Tool: Lemma \ref{bnj08}]
$\bullet$ First we show that it is sufficient to prove the
existence of an integer $c\in [-N,N]\setminus\{0\}$ instead of
$c\in [N].$ 
Observe that it is true if there exists a subset
$T\subset S\subset [N]_0^r$ with $|T|\ge \delta_r N^r$ such that
$T$ is symmetric with respect to some $x_T=x\in ({1\over
2}[2N]_0)^r:=\{{z\over 2}| z\in [2N]_0\}^r$ (i.e., for any $t\in
T$ there is a $t'\in T$ with ${1\over 2}(t+t')=x$) where
$\delta_r>0$ is a constant independent of $N.$ 
Randomly picking
a point $x\in ({1\over 2}[2N]_0)^r,$ the expected number of pairs
$s,s'\in S$ with $s+s'=2x$ is ${|S|\choose 2}/(2N)^r \ge
0.49{\delta^2}N^{2r}/(2N)^r={0.49\over 2^r} \delta^2 N^r.$ Thus
there exists the desired $T$ with $|T|\ge {0.98\over 2^r}\delta^2
N^r.$
\par
$\bullet$ By the above remark, it easily follows from
\underline{Lemma \ref{bnj08}} that the theorem holds when
$F\subset B_r:=E_r\cup\{0=(0,\cdots,0)\}$, by ignoring the
$0^{th}$ coordinate.
\par
$\bullet$ Let $\delta,\, r,\, F,$ and $S$ be given as in the
theorem. Without loss of generality, $F$ can be written as
$[m]_0^r=\{(x_1,\cdots,x_r)| x_i\in [m]_0\}$ for a constant $m.$
Let $r':=|F|-1=m^r-1.$ 
Take a linear map
$\phi:\realset^{r'}\to\realset^r$ such that the restriction
$\phi|_{B_{r'}}$ is a bijection from $B_{r'} $ to $F$ with
$\phi(0)=0.$ 
Define $S'\subset [N]_0^{r'}$ by
$S':=\phi^{-1}(S)\cap [N]_0^{r'}=\{z \,|\, \phi(z)\in S\}.$ 
Clearly $\phi^{-1}(s)$ forms an $(r'-r)$-dimensional 
linear
subspace of $ \realset^{r'}$ for any $s\in S$, by observing the rank of 
an $r\times r'$-matrix.
Then it is straightforward to
see that there exists a constant $\delta'=\delta'(\delta,r,m)$
such that $|S'|\ge \delta'N^{r'}$. 
Taking $N$ large, the last
paragraph yields
 $a\in [N]_0^{r'}$ and $c>0$ such that
$a+cB_{r'}\subset S'.$
Thus $S\supset\phi(a+cB_{r'})=\phi(a)+c\phi(B_{r'})=
\phi(a)+cF,$ completing the proof.
\qed
\section{Remarks}
Let $F$ be a $k$-uniform  ($2$-colored: black and invisible)
hypergraph. Denote by ${\rm ex}^{(k)}(n,F)$ the maximum number of
black edges of a $k$-uniform ($2$-colored: black and white)
hypergraph on exactly $n$ vertices
 with no copy of $F$ as a subgraph. By an easy modification of the proof of our removal lemma,
we can easily show a hypergraph version of the Erd\H{o}s-Stone
theorem.
\begin{prop}[A hypergraph version of the Erd\H{o}s-Stone theorem]
Let $F, F_0$ be any $k$-uniform hypergraphs such that
 $F$ is a \lq blow-up\rq\ of $F_0$ (i.e.,
there exists a map from the vertex set $V(F)$ to $V(F_0)$ such
that
 each (black) edge of $F$ is
mapped to a (black) edge of $F_0$). Then ${\rm ex}^{(k)}(n, F)\le
{\rm ex}^{(k)}(n,F_0)+o(n^k).$
\end{prop}
R\"{o}dl and Skokan \cite{RSk06} have already shown the above for
black-only $F_0$ (i.e., $F_0=K_r^{(k)}$) by adding extra arguments
 to a removal lemma.
Although it should not be hard to obtain the above
by previously known techniques, ours is a direct and shorter proof.
\par
It is worthwhile to note that not only the way of regularizing  but also
the construction of
the error function (\ref{lj25}) is quite simple and clear in our proof.
It is easy to find
a simple $O(1)$-time random algorithm by which we can approximately
grasp the entire hypergraph ${\bf G}.$
\par
Alon et al. \cite{AFNS} discussed the relation between Regularity
Lemma and Property Testing for ordinary graphs. Although their
proof is conceptually clear, many of their technical details may
come from their problem setting (non-partiteness). In order to
understand the essential relation between Regularization
(Regularity Lemma) and Property Testing, it may be even easier and
more natural to consider them on partite hypergraphs rather than
on nonpartite ordinary graphs. {\em Property Testing and
Regularization are essentially equivalent.} They are all about
random samplings. If there exists a difference between the two, it
is
\begin{em}
whether the number of random vertex samplings
is (PT) a fixed constant or (R) bounded by a constant but chosen randomly.
\end{em}
The above difference is essentially insignificant, as far as we do
not consider the sizes of constants. Property Testing is stronger
than Regularization in the sense that a (non-canonical) property
tester can ignore some random number of vertex samples after
choosing the vertices.\footnote{ Canonical property testing
chooses (a fixed number of) vertices at random, but once the
vertices are chosen, it outputs its answer deterministically.
Therefore, at first sight, canonical property testing may be
weaker. However as seen in \cite{AFKS}\cite[Th.2]{GT}, for any
given non-canonical property test, there exists a canonical
property test which is equivalent to it. (Its derandomizing
process is easy, since the sampling size of a non-canonical tester
is a constant. The canonical tester repeats the samplings many
(but a constant number of) times.
 Then it computes the probability that the noncanonical tester accepts
for each sampling. The canonical tester accepts iff the average 
of the
probabilities is at least $1/2.$ )} On the other hand,
Regularization is stronger than Property Testing in the sense that
Regularization \lq knows\rq\ the number of copies of all
fixed-sized subgraphs approximately. (If there is another
difference, the Property Tester outputs one of only two choices
(YES/NO), while Regularization can output some of a constant
number of choices; also see \cite{I06m}).
\par
Therefore our result on hypergraph regularization is not a simple
extension of graph regularization. It helps our understanding of
regularization (and property testing)  both for graphs and
hypergraphs.

\end{document}